\documentclass[]{article}
\usepackage[a4paper,margin=3cm]{geometry}

\usepackage[utf8]{inputenc}
\usepackage[T1]{fontenc}
\usepackage{lmodern}
\usepackage[ngerman,english]{babel}
\usepackage{amsmath,amsfonts,amssymb,mathtools}
\usepackage{import}
\usepackage{authblk}
\usepackage[hidelinks]{hyperref}
\usepackage[noabbrev,nameinlink]{cleveref}
\crefname{equation}{}{}
\usepackage{graphicx}
\usepackage{tabularx}

\usepackage{float,subcaption}
\usepackage{tikz}
\usepackage{algorithm}
\usepackage{algorithmicx}
\usepackage{algpseudocode}
\usepackage{booktabs}
\usetikzlibrary{decorations.pathreplacing,spy}

\renewcommand\vec{\boldsymbol}

\DeclareMathOperator{\Sym}{Sym}
\DeclareMathOperator{\Div}{div}

\DeclareMathOperator{\Tr}{Tr}

\newcommand{\R}{\mathbb{R}}
\newcommand{\F}{F}

\newcommand{\Lagr}{{\mathcal{L}}}
\newcommand{\Gobs}{{\Gamma_{\text{obs}}}}
\newcommand{\Gout}{{\Gamma_{\text{out}}}}
\newcommand{\Gin}{{\Gamma_{\text{in}}}}
\newcommand{\Gwall}{{\Gamma_{\text{wall}}}}
\newcommand{\Oobs}{{\Omega_{\text{obs}}}}

\newcommand{\holdall}{{G}}
\newcommand{\bc}{{\text{bc}}}
\newcommand{\vol}{{\text{vol}}}

\newcommand{\test}[1]{{\delta_{#1}}}
\newcommand{\mult}[1]{{\lambda_{#1}}}
\newcommand{\vectest}[1]{{\vec{\delta}_{#1}}}
\newcommand{\vecmult}[1]{{\vec{\lambda}_{#1}}}

\newcommand{\control}{{u}}
\newcommand{\edet}{{\eta_{\text{det}}}}

\newcommand{\vel}{\vec{v}}
\newcommand{\defor}{\vec{w}}
\newcommand{\Transform}[1]{(D#1\InvF)}
\newcommand{\InvF}{(D\F)^{-1}}
\newcommand{\Det}{\det(DF)}

\newcommand{\Om}{\Omega}

\newcommand{\adjdef}{\mult{\defor}}
\newcommand{\adjvel}{\mult{\vel}}
\newcommand{\adjpress}{\mult{p}}
\newcommand{\adjbc}{\mult{\bc}}
\newcommand{\adjvol}{\mult{\vol}}
\newcommand{\etaub}{{\eta_\mathrm{ub}}}
\newcommand{\etalb}{{\eta_\mathrm{lb}}}

\title{Fluid dynamic shape optimization using self-adapting nonlinear extension operators with multigrid preconditioners}

\author[1]{Jose Pinzon}
\author[1]{Martin Siebenborn}

\affil[1]{{Department of Mathematics, University Hamburg, Bundesstr. 55, 20146 Hamburg, Germany}}

\begin{document}
	\maketitle
	
	\begin{abstract}
		In this article we propose a scalable shape optimization algorithm which is tailored for large scale problems and geometries represented by hierarchically refined meshes.
		Weak scalability and grid independent convergence is achieved via a combination of multigrid schemes for the simulation of the PDEs and quasi Newton methods on the optimization side.
		For this purpose a self-adapting, nonlinear extension operator is proposed within the framework of the method of mappings.
		This operator is demonstrated to identify critical regions in the reference configuration where geometric singularities have to arise or vanish.
		Thereby the set of admissible transformations is adapted to the underlying shape optimization situation.
		The performance of the proposed method is demonstrated for the example of drag minimization of an obstacle within a stationary, incompressible Navier-Stokes flow.
	\end{abstract}
	

\par
\section{Introduction}\label{sec:intro}
%
Shape optimization is a mathematical tool to obtain an optimal contour for a randomly-shaped obstacle subject to physical phenomena described by a partial differential equation (PDE). 
In the classical sense, this is achieved by evaluation of sensitivities of a shape functional $j(y,\Om)$ for a set of admissible shapes $\Omega \in G_\text{adm}$. 
The functional is constrained by one or several PDEs, among them a state equation $E(y,\Om)=0$, which is fulfilled by the state variable $y$. 
In this article, we focus on incompressible flow described by Navier-Stokes equations where the objective is to minimize the energy dissipation around an obstacle in terms of its shape and additional geometrical constraints. 
Building on a long history ranging back several decades (see for instance \cite{jameson2003aerodynamic,giles2000introduction,mohammadi2010applied}), the field of shape optimization with fluid dynamics applications is still very active today following various approaches (e.g.\ \cite{schmidt2013three,mueller2021novel,garcke2016stable,FLUU16}).
For an overview of shape optimization constrained by PDEs we refer the reader to \cite{sokolowski2012introduction,allaire2021shape,delfour2001shapes}.

The iterative application of deformation fields to a finite element mesh can lead to distortions and loss of mesh quality, as studied by~\cite{dokken2019shape,EHLG18}. 
This becomes particularly disruptive for numerical algorithms if there are large deformations leading from reference to optimal configuration.
Several approaches, especially in recent studies, have been proposed to prevent this. 
For instance, the use of penalized deformation gradients in interface-separated domains helps maintain mesh quality but might still lead to element overlaps when taken to the limit \cite{VogelSiebenborn2021}.
Other approaches rely on remeshing the domain, as for instance in \cite{Wilke2005}. 
More recent efforts on this area make use of pre-shape calculus to allow for the simultaneous optimization of both the shape and mesh quality of the grid \cite{preshape_a2021}.

Although the variable of the optimization problem is only the contour of the shape, the surrounding space plays a crucial role since it describes the domain for the physical effects.
Due to the Hadamard-Zolésio structure theorem (see for instance \cite{sokolowski2012introduction}) changes of the objective function under consideration can be traced back to deformations of the shape, which are orthogonal to its boundary.
This has been recognized as a source for decreasing mesh qualities and is addressed by many authors, for instance by also allowing displacements tangent to the shape surface (cf. \cite{preshape_b2021}).
In contrast to the statement of the structure theorem, from a computational point of view, it can be favorable to extend surface deformations into the entire surrounding domain instead of building a new discretization around the improved shape after each descent step, i.e. avoid remeshing on each new iteration.
In recent works it has become popular to reuse the domain around the shape, which describes the domain of the PDE, for the representation of Sobolev gradients (e.g. \cite{schulz2016efficient}).
Typically, elliptic equations are solved in this domain in order to represent the shape sensitivity as a gradient with respect to a Hilbert space defined therein \cite{DokkenFunke20,harbrecht2013numerical,langer2015shape}.
The benefit of this approach is that the resulting deformation field not only serves as a deformation to the obstacle boundary, but can also be utilized as a domain deformation.
Thus, the discretization for the next optimization step is obtained without additional computational cost.

At this point, two different approaches can be distinguished.
On the one hand, the computed gradient can be used directly as a deformation to the domain after each optimization step which can be seen as changing the reference domain iteratively.
On the other, in the so called method of mappings \cite{MuSi}, the reference domain is fixed and the shape optimization problem is interpreted via the theory of optimal control.
These are implemented through the definition of a variable around the surface of the shape to be optimized and its connection with the deformation field affecting the whole domain through a so called extension operator. 
The solution of which results in the optimal deformation field for both the target shape and its surrounding domain.
For an application of this method to aerodynamic shape optimization see for instance \cite{FLUU16}.

We can oppose these two approaches as
\begin{equation}
	\begin{aligned}
	&\min\limits_{\Omega \in \holdall_\text{adm}} && j(y, \Omega) &\qquad\qquad& \min\limits_{\F \in \F_\text{adm}} && j(y, \F(\Omega))  \\
	&\quad \text{s.t.} && e(y, \Omega) = 0 &\qquad\qquad& \quad\text{s.t.}\quad && e(y, \F(\Omega)) = 0
	\end{aligned}
\end{equation}
where either a set of admissible domains $\holdall_\text{adm}$ or a set of admissible transformations $\F_\text{adm}$ has to be defined.
A link between these two can then be established via
\begin{equation}
	\holdall_\text{adm} := \left\{ \F(\Omega):\, \F \in \F_\text{adm}\right\}
\end{equation}
in terms of a given reference configuration $\Omega$.\\

The approach we propose here is based on the research carried out in \cite{haubner2021continuous}. 
Compared to iteratively updating the shape, it offers the possibility to require properties of the deformation from reference to optimal shape.
Moreover, it reformulates the optimization over a set of admissible transformations $F \in F_\text{adm}$, which enables us to carry out the optimization procedure on the reference domain. 
Additionally, it has been documented that such extension operators are possible without the need for heuristic parameter tuning.
In \cite{onyshkevych2021} an additional nonlinear term is introduced to the elliptic extension equation allowing for large deformations while preserving mesh quality, and preventing element self-intersections and degeneration. 
In shape optimization this occurs when trying to obtain an optimal shape for an obstacle surrounded by media, for which the creation or removal of a singularity on the obstacle's boundary is necessary.

In the present work we focus on applying parallel solvers for the solution of PDEs in large distributed-memory systems.
This stems from the fact that the discretizations will lead to a very large number of degrees of freedom (DoFs), for which the application of the geometric multigrid method (see for instance \cite{Hackbusch85}) guarantees mesh-independent convergence on the simulations.
Its application within the context of parallel computing towards the solution of PDEs is an area of ongoing research \cite{Reiter2014gmg,gmeiner2014parallel,Baker11}.
The feature of mesh-independent convergence is a necessary condition towards weak scalability of the entire optimization algorithm, which is why in this article we apply the multigrid method as a preconditioner for the solution of the PDE constraints. 
This requires to provide a sequence of hierarchically refined discretizations.
However, the shape optimization problem is a fine grid problem, which means that the contour of the obstacle has to be representable within the entire grid hierarchy.
This leads to undesired effects and conceptual challenges that have been addressed for instance in \cite{naegel2015scalable,siebenborn2017algorithmic,pinzon2020parallel}.
The finest grid in the hierarchy thus typically represents a high-resolution discretization of a polygonal shape, which -- besides the aforementioned sources for mesh degeneracy -- also introduces challenges to the shape optimization.
This is due to the fact that it results in singularities (e.g.\ edges and corners) that are visible to the discretization.
Particularly for fluid dynamics applications this is problematic.
If the regularity of the considered domain transformations, i.e.\ the descent steps of the optimization, is then too high, it is not possible to achieve smooth optimal shapes.
Related to this is the problem, that, away from the hierarchical grid structure, if one starts with a smooth reference shape, it is not possible to have singularities in the optimal configuration.

In this work, we propose an approach that is able to identify these regions and adapt the set of admissible transformations $\F_\text{adm}$ as part of the optimization problem.
In \cite{onyshkevych2021} it is illustrated how adding a nonlinear convection term to the extension model that defines $\F_\text{adm}$ affects forming singularities in optimal shapes.
We thus study in this article how the non-linearity can be adjusted according to the shape of the reference domain.

The rest of this article is structured as follows: In \cref{sec:background} the optimization problem is formulated and the underlying fluid experiment is outlined.
\Cref{sec:algorithm} is devoted to the optimization algorithm and the computation of the reduced gradient via the adjoint method.
Subsequently, in \cref{sec:results} the performance of the proposed method is discussed by presenting numerical tests.
The article closes with a conclusion in \cref{sec:conclusion}.

\section{Problem formulation and mathematical background}\label{sec:background}
\begin{figure}[tbp]
	\begin{center}
		\scalebox{0.8}{\includegraphics{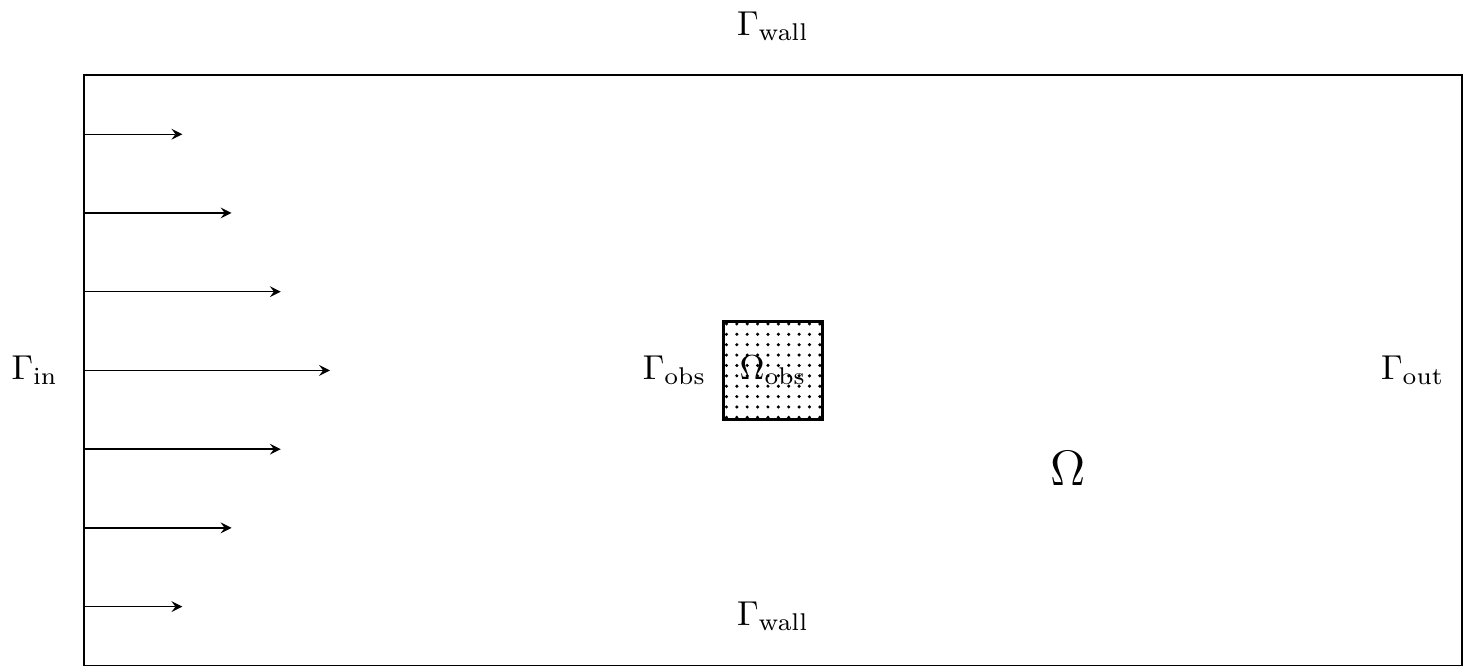}}
	\end{center}
	\caption{2d holdall reference domain of the flow field with square obstacle.}
	\label{fig:2d_holdall}
\end{figure}

The model problem under consideration is sketched in \cref{fig:2d_holdall} in a bounded holdall domain $\holdall := \Omega \cup \Oobs$, where $\Omega$ is assumed to have a Lipschitz boundary. 
In $\Omega$ we consider a stationary, incompressible flow.
It surrounds an obstacle $\Oobs$ with variable boundary $\Gobs$, but fixed volume and barycentric coordinates. 
Throughout this article, the original setting of the domain will be referred to as the reference configuration or domain.

In this section we present the theoretical background that culminates with the algorithm presented in \cref{sec:algorithm}. 
The problem is first formulated in terms of classical shape optimization, to be then reformulated as an optimal control problem. 
Later on, it is pulled back to the reference domain through the method of mappings. 
The extension operator used in our approach is described through its weak formulation to be able to formulate the augmented Lagrange approach later used. 
Finally, the Lagrangian is used to obtain the sensitivities necessary for the descent direction and for the approximation to the Hessian used. 
For an in-depth discussion of the underlying theory we refer the reader, as previously mentioned, to \cite{onyshkevych2021,BaLiUl,haubner2021continuous}, the use of the adjoint method in fluid dynamics can be reviewed in \cite{giles2000introduction,jameson2003aerodynamic}, and \cite{hinzeulbrich2009} can be consulted for a mathematical review of the theory here used.\\
Let $X = L^2(\Gobs) \times L^2(\Omega)$, $0 < \etalb < \etaub$, $b > 0, \alpha > 0, \theta > 0$ and consider the optimization problem
\begin{equation}
	\min\limits_{(\control, \eta) \in X}\quad j(y, \F(\Omega)) + \tfrac{\alpha}{2} \| \control \|_{L^2(\Gobs)}^2 + \tfrac{\theta}{2}\| \eta - \tfrac{1}{2}(\etaub+\etalb) \|_{L^2(\Omega)}^2\label{eq:opt_sys_0}
\end{equation}
\begin{align}
	\text{s.t.}\quad& e(y, \F(\Omega)) = 0\label{eq:opt_sys_1}\\
	& \F = \mathrm{id} + \defor\label{eq:opt_sys_2}\\
	& \defor = S(\eta, \control, \Omega)\label{eq:opt_sys_3}\\
	& \det(DF) \geq b   && \text{in}\; \Omega\label{eq:opt_sys_4}\\
	& \etalb \leq \eta \leq \etaub && \text{in}\; \Omega\label{eq:opt_sys_5}\\
	& g(\defor) = 0.\label{eq:opt_sys_6}
\end{align}
where $g(\defor)$ represents geometric constraints.
$S$ denotes a so-called extension operator, which links the boundary control variable $\control \in L^2(\Gobs)$ to a displacement field $\defor: \Omega \to \R^d$.
Examples of possible choices of $S$ are given and investigated in \cite{haubner2021continuous}.
Particularly, the resulting regularity of the domain transformation $\F$ is discussed therein.
Here we enrich this operator with an additional control variable $\eta \in L^2(\Omega)$, which plays the role of a nonlinearity switch.
In the following we assume that $S$ is defined such that $\defor|_{\Gin\cup\Gwall\cup\Gout} = 0$ a.e.

In the experiment presented here, the geometric constraints require the barycenter and volume to be the origin and constant, respectively.
This excludes trivial solutions where the obstacle shrinks to a point or moves towards a position where the objective functional is minimized.
Thus, the principal geometric constraints are given by 
\begin{equation}
	\hat{g}(\defor) = \left( \int_{F(\Omega)} 1 \, dx - \int_\Omega 1 \, dx , \frac{\int_{F(\Omega)} x\, dx}{\int_{F(\Omega)} 1 \, dx} - \frac{\int_\Omega x\, dx}{\int_\Omega 1 \, dx}\right)^\top
	\label{eq:pre_geometric_constraints}
\end{equation}
This simplifies to
\begin{equation}
	g(\defor) = \left(\int_\Omega \det(DF) - 1 \, dx , \int_\Omega F(x) \det(DF) \, dx\right)^\top,
	\label{eq:geometric_constraints}
\end{equation}
assuming, without loss of generality, that the barycenter of the reference domain is $0\in\R^d$ and the volume is precisely retained.

The condition \eqref{eq:opt_sys_4} is approximated via a non-smooth penalty term.
This results in an approximation via the objective function
\begin{multline}
	J(y, \control, \defor, \eta) = j(y, \F(\Omega)) + \tfrac{\alpha}{2} \| \control \|_{L^2(\Gobs)}^2 + \tfrac{\theta}{2}\| \eta - \tfrac{1}{2}(\etaub+\etalb) \|_{L^2(\Omega)}^2\\
	 + \tfrac{\beta}{2}\|(b - \det(DF))^+\|^2_{L^2(\Omega)}.
\label{eq:objective_func}
\end{multline}

In contrast to the PDE constraints \cref{eq:opt_sys_2,eq:opt_sys_3,eq:opt_sys_4}, the geometric constraints \cref{eq:opt_sys_6} are fixed dimensional (here it is $d+1$ where $d\in \lbrace2, 3\rbrace$).
Thus, the multipliers associated to these conditions are not variables in the finite element space but a $d+1$-dimensional vector.
This is incorporated into the optimization algorithm in form of an augmented Lagrange approach.
This leads to the augmented objective function
\begin{equation}
	J_\text{aug}(y, \control, \defor, \eta) := J(y, \control, \defor, \eta) + \tau \|g(\defor)\|_2^2
\label{eq:augmented_objective_func}
\end{equation}
where $\tau > 0$ is a penalty factor for the geometric constraints.
The basic concept of the augmented Lagrange method is to optimize the objective \cref{eq:augmented_objective_func}.
By contrast to a pure penalty method, the geometric constraints \cref{eq:opt_sys_6} are not entirely moved to the objective function, but the corresponding multipliers $\vecmult{g}$ are assumed to be approximately known and iteratively updated.

We consider the PDE constraint $e(y,\F(\Omega))$ to be the stationary, incompressible Navier-Stokes equations in terms of velocity and pressure $(\vel,p)$.
In the following it is distinguished between PDE solutions defined on the reference domain $\Omega$ denoted by $(\vel,p)$ and on the transformed domain $\F(\Omega)$ as $(\hat{\vel},\hat{p})$. 
We thus consider
\begin{align}
	-\nu \Delta \hat{\vel} + (\hat{\vel}\cdot \nabla) \hat{\vel} + \nabla \hat{p} = 0 &\qquad \text{in} \; \F(\Omega)\label{eq:ns_1}\\
	\Div \hat{\vel} = 0 &\qquad \text{in} \; \F(\Omega)\label{eq:ns_2}\\
	\hat{\vel} = \vel_\infty &\qquad \text{on} \; \Gin\label{eq:ns_3}\\
	\hat{\vel} = 0 &\qquad \text{on} \; \Gobs \cup \Gwall\label{eq:ns_4}\\
	\hat{p}n - \nu \frac{\partial \hat{\vel}}{\partial n} = 0  &\qquad \text{on} \; \Gout,\label{eq:ns_5}\\ 
\end{align}
where for compatibility it is assumed that $\int_\Gin \vel_\infty\cdot n\, ds=0$ holds for the inflow velocity profile $\vel_\infty$.
Notice that the boundaries $\Gin, \Gout, \Gwall$ in \cref{eq:ns_2,eq:ns_3,eq:ns_4} are unchanged since the displacement $\defor$ is assumed to vanish here.
This assumption reflects that the outer boundaries of the experiment domain are not a variable in the optimization problem.

The variational formulation of the PDE \crefrange{eq:ns_1}{eq:ns_5} pulled back to the reference domain $\Omega$ is given by: Find $(\vel, p)\in V \times Q$ such that for all $(\vectest{\vel}, \test{p})\in V_0 \times Q$ it holds
\begin{align}
	\int_\Omega \nu (D\vel \InvF) : (D\vectest{\vel} \InvF) + (D\vel \InvF\vel)\cdot \vectest{\vel} - p \Tr(D\vectest{\vel}\InvF)\Det\, dx = 0\label{eq:ns_variational_1}\\
	-\int_\Omega \test{p} \Tr(D\vel\InvF)\Det\, dx=0\label{eq:ns_variational_2},
\end{align}
where trial and test functions are chosen in
\begin{equation}
	\begin{aligned}
		V &:= \left\{ \vel \in H^1(\Omega, \R^d): v|_\Gin = \vel_\infty, v|_{\Gwall \cup \Gobs} = 0 \text{ a.e.}\right\},\\
		V_0 &:= V \;\text{ with }\; \vel_\infty = 0, \\
		Q &:= \left\{ p \in L^2(\Omega): \int_\Omega p\, dx = 0 \right\}.
	\end{aligned}
\end{equation}

In the experiment considered in this work the physical part of the objective function \cref{eq:objective_func} is given by the energy dissipation in terms of the velocity $\vel$, thus $y=(\vel, p)$ and
\begin{equation}
	\hat{j}(\hat{\vel}, \F(\Omega)) = \nu \int_{\F(\Omega)} D\hat{\vel} : D\hat{\vel}\, dx,
\end{equation}
which can be pulled back to the reference domain $\Omega$ as
\begin{equation}
	j(\vel, \defor) = \nu \int_{\Omega} (D\vel\InvF) : (D\vel\InvF)\,\Det dx.
\end{equation}

The extension $S(\eta, \control, \Omega)$ is defined to be the solution operator of the PDE
\begin{align}
	\Div (D \defor + D \defor^\top) + \eta(\defor\cdot \nabla)\defor = 0 &\qquad \text{in} \; \Omega\label{eq:extension_1}\\
	(D \defor + D \defor^\top)n = \control \vec{n} &\qquad \text{on} \; \Oobs\label{eq:extension_2.1}\\
	\defor = 0 &\qquad \text{on} \; \Gin \cup \Gout \cup \Gwall.\label{eq:extension_2.2}
\end{align}
Consider the space 
\begin{equation}
	W := \left\{ \defor \in H^1(\Omega, \R^d): \defor|_{\partial \Omega \setminus \Gobs} = 0 \text{ a.e.} \right\}.
\end{equation}
Then the variational formulation of \cref{eq:extension_1,eq:extension_2.1,eq:extension_2.2} is obtained by:
Find $\defor \in W$ such that for all $\vectest{\defor} \in W$ it holds
\begin{equation}
	\int_\Omega \Sym(D\defor):D\vectest{\defor} + \eta (D\defor\, \defor) \cdot \vectest{\defor}\, dx = \int_\Gobs u \vec{n} \cdot \vectest{\defor}\,ds.
	\label{eq:weak_deform_equation}
\end{equation}

Finally, we can formulate the approximate optimization problem, which is then solved via the augmented Lagrange approach in \cref{sec:algorithm}, as
\begin{align}
	\min\limits_{(\control, \eta) \in X}&\quad J_\text{aug}(y, \control, \defor, \eta) &&\label{eq:opt_sys_final_0} \\
	\text{s.t.}\quad& \text{\cref{eq:ns_1,eq:ns_2,eq:ns_3,eq:ns_4,eq:ns_5,eq:weak_deform_equation}} &&\label{eq:opt_sys_final_0.1} \\
	& \F = \mathrm{id} + \defor && \label{eq:opt_sys_final_1}\\
	& \etalb \leq \eta \leq \etaub && \text{in}\; \Omega\label{eq:opt_sys_final_2}\\
	& g(\defor) = 0, && \label{eq:opt_sys_final_3}
\end{align}
where the multipliers for condtions \cref{eq:opt_sys_final_3} are assumed to be known in each iteration.

In order to formulate a gradient-based decent algorithm, we have to compute sensitivities of the final objective function $J_\text{aug}$ in \cref{eq:augmented_objective_func} with respect to the variables $(\control, \eta)$.
This means to differentiate the chain of mappings
\begin{equation}
	(\control, \eta) \mapsto \defor \mapsto (\vel, p) \mapsto J_\text{aug}(\vel , \control, \defor, \eta)
\end{equation}
and obtain the sensitivities in reverse order
\begin{equation}
	J_\text{aug}(\vel , \control, \defor, \eta) \mapsto (\vecmult{\vel}, \mult{p}) \mapsto \vecmult{\defor} \mapsto (\mult{u}, \mult{\eta}).
\end{equation}
Access to the adjoint gradient formulation can be obtained via the corresponding Lagrangian, which is given by
\begin{multline}
		\Lagr (\vec{\defor},\vec{v},p,\control,\eta,\adjdef,\adjvel,\adjpress,\adjbc\adjvol) = \frac{\nu}{2}\int_{\Om}\Transform{\vel}:\Transform{\vel}\Det \,dx\\
		+ \frac{\alpha}{2}\int_{\Gobs}\control^2 \,ds + \frac{\beta}{2}\int_{\Om}((\edet-\Det)^+)^2 \,dx + \frac{\theta}{2}\int_\Omega \left(\eta - \tfrac{1}{2}(\etaub+\etalb)\right)^2 \, dx + \tau \|g(\defor)\|_2^2 \\
		+ \int_{\Om}[\nu\Transform{\vel}:\Transform{\mult{\vel}} + \Transform{\vel}\cdot \mult{\vel} - p \Tr\Transform{\mult{\vel}} ]\Det \,dx\\
		- \int_{\Om}\mult{p}\Tr\Transform{\vel}\Det \,dx + \int_{\Om}[(D\defor + D\defor^\top):D\mult{\defor} + \eta(D\defor\cdot\defor)]\,dx\\ 
		- \int_{\Gobs}\control\vec{n}\cdot\mult{\defor}\,ds + \mult{\bc}\cdot\int_{\Om}\F(x)\Det \,dx + \mult{\vol}\int_{\Om}(\Det-1)\,dx
	\label{eq:lagrangian}
\end{multline}
under the assumption that the barycenter of $\Om$ is $0 \in \R^d$.

From the Lagrangian \cref{eq:lagrangian} the adjoint Navier-Stokes equations follow as: Find $(\vecmult{\vel}, \mult{p})\in V_0 \times Q$ such that for all $(\vectest{\vel}, \test{p})\in V_0 \times Q$ it holds

\begin{multline}
	\int_\Omega \nu (D\vecmult{\vel} \InvF) : (D\vectest{\vel} \InvF) + (D\vectest{\vel} \InvF\vel)\cdot \vecmult{\vel}\\
	 + (D\vel \InvF\vectest{\vel})\cdot \vecmult{\vel} - \mult{p} \Tr(D\vectest{\vel}\InvF)\Det\, dx = 0\label{eq:ns_adjoint_variational_1}
\end{multline}
\begin{equation}
		-\int_\Omega \test{p} \Tr(D\vecmult{\vel}\InvF)\Det\, dx=0\label{eq:ns_adjoint_variational_2},
\end{equation}
The adjoint displacement equation is obtained by: Find $\vecmult{\defor} \in W$ such that for all $\vectest{\adjdef} \in W$ it holds
\begin{equation}
	\int_\Omega \Sym(D\vecmult{\defor}):D\vectest{\adjdef} + \eta (D\defor\, \vecmult{\defor}) \cdot \vectest{\adjdef}\, dx = R(\defor, \vel, p, \vecmult{\vel}, \mult{p}).
	\label{eq:weak_adjoint_deform_equation}
\end{equation}
In \cref{eq:weak_adjoint_deform_equation} $R$ denotes the derivative of the Lagrangian \cref{eq:lagrangian} w.r.t.\ $\defor$.
This is obtained after straightforward computations and omitted here for the sake of brevity.
Finally, the reduced gradient is obtained as: Find $(\gamma, \kappa) \in X$ such that for all $(\test{\control}, \test{\eta})\in X$ it holds
\begin{align}
	\int_\Gobs \gamma \test{\gamma}  +\alpha \control \test{\control} - \vecmult{\defor} \cdot n \test{\control} \, ds &= 0,
	\label{eq:control_gradient}\\
	\int_\Omega \kappa \test{\eta}  +\theta \left(\eta - \tfrac{1}{2}(\etaub+\etalb)\right) \test{\eta} - (D \defor \cdot \defor) \cdot \vecmult{\defor} \test{\eta} \, ds &= 0.
	\label{eq:extension_gradient}
\end{align}
With the sensitivity equations \cref{eq:control_gradient,eq:extension_gradient} we are now prepared to apply a descent method.
\newline
\section{Optimization Algorithm}\label{sec:algorithm}
In \cref{sec:background} we present the approximate optimization problem \cref{eq:opt_sys_final_0,eq:opt_sys_final_0.1,eq:opt_sys_final_1,eq:opt_sys_final_2,eq:opt_sys_final_3}, which is solved via the augmented Lagrange approach shown in \cref{alg::aug_lag_alg}. 
An initial guess is given to the Lagrange multipliers $\vecmult{g}$ associated to the geometrical constraints. 
These in turn are iteratively updated in each optimization step subject to the condition that the norm of the defect the geometrical constraints is smaller than a prescribed tolerance $\epsilon_g > 0$.

Most of the computational time is consumed for solving the PDE systems presented in \cref{sec:background}. 
This is carried out by \cref{alg::gradient} in a block-wise manner, where the output consists of the new displacement field to update the transformation \cref{eq:opt_sys_2}, as well as the results of the reduced gradient $\nabla \control^{k,\ell}, \nabla\eta^{k,\ell}$, which will be further used to obtain updates for the current control and extension factor, $\control^{k,\ell+1}, \eta^{k,\ell+1}$, respectively. 
A reference for the reduced gradient method can be found for instance in \cite{hinzeulbrich2009}.
\begin{algorithm}[!ht]
\caption{Augmented Lagrange Outer Optimization Algorithm}
\label{alg::aug_lag_alg}
	\begin{algorithmic}[1]
	\Require $\vecmult{g}, \mult{\mathrm{inc}} > 0, 0 \leq \etalb \leq \eta \leq \etaub,  b > 0, \epsilon_g > 0, \tau > 0, \tau_\mathrm{inc} > 0,\epsilon_\mathrm{outer} > 0, \epsilon_\mathrm{inner} > 0$
	\Repeat
	\State $\ell \gets 0$
	\Repeat
	\State $(w, \nabla \control^{k,\ell}, \nabla\eta^{k,\ell}) \gets$ \Call{Reduced Gradients}{$g_\mathrm{def}, \control^{k,\ell}, \eta^{k,\ell}$}
	\State $(\control^{k,\ell+1}, \eta^{k,\ell+1}) \gets$  \Call{lBFGS-B}{$m, \ell, g_\mathrm{def}, \tau, \control^{k,\ell}, \eta^{k,\ell}, \nabla \control^{k,\ell}, \nabla\eta^{k,\ell}$}
	\State $g_\mathrm{def} \gets g(w) - g(0)$
	\State $\ell \gets \ell + 1$
	\Until{$\left\|\left(\nabla \control^{k,\ell-1}, P_{(\etalb, \etaub)}\left( \eta^{k,\ell-1} - \nabla \eta^{k,\ell-1}\right) - \eta^{k,\ell-1}\right) \right\|_X < \epsilon_\mathrm{inner}$}
	\If{$\|g_\mathrm{def}\|_2 < \epsilon_g$}
	\State $\tau \gets \tau_\mathrm{inc} \tau$
	\Else
	\State $\vecmult{g} \gets \vecmult{g} + \mult{\mathrm{inc}} g_\mathrm{def}$
	\EndIf
	\State $k \gets k+1$
	\Until{$k\geq1$ \;\textbf{and}\; $\|\control^{k,\ell} - \control^{k-1,\ell}\|_{L^2(\Gobs)} < \epsilon_\mathrm{outer}$}
	\end{algorithmic}
\end{algorithm}

The objective function is non-differentiable due to the presence of the positive part mapping in $R$ in~\cref{eq:weak_adjoint_deform_equation}, which is discussed in depth in~\cite{haubner2021continuous}. Moreover, a discussion on quasi-Newton methods for semi-smooth objective functions and  can be found in \cite{mannel2020hybrid}.

\begin{algorithm}[!ht]
	\caption{Limited Memory BFGS Algorithm with Box Constraints Inner Iteration}
	\label{alg::lBFGS}
	\begin{algorithmic}[1]
		\Function{lBFGS-B}{$m, \ell, g_\mathrm{def}, \tau, w, \control^\ell, \eta^\ell, \nabla \control^\ell, \nabla \eta^\ell$}
		\State Determine $\chi_\eta$ according to \cref{eq:eta_indicator} and $(\cdot, \cdot)_{\hat{X}}$ according to \cref{eq:modified_inner_prod}
		\State $q = (q_1,q_2) \gets (\nabla \control^\ell, \nabla \eta^\ell)$
		\If{$\ell > 0$}
		\State $s^{\ell-1} \gets (\control^\ell, \eta^\ell) - (\control^{\ell-1}, \eta^{\ell-1})$
		\State $z^{\ell-1} \gets (\nabla \control^\ell, \nabla \eta^\ell) - (\nabla \control^{\ell-1}, \nabla \eta^{\ell-1})$
		\For{$i=\ell-1,\dots, \max\{\ell-m, 0\}$}
			\State $\alpha_i \gets \rho_i (s^i, q)_{\hat{X}}$
			\State $q \gets q - \alpha_i z^i$
		\EndFor
		\State $q \gets \frac{(s^{\ell-1}, z^{\ell-1})_{\hat{X}}}{(z^{\ell-1}, z^{\ell-1})_{\hat{X}}} q$
		\For{$i=\max\{\ell-m, 0\},\dots, \ell-1$}
			\State $q \gets q + (\alpha_i - \rho_i (y^i, q)_{\hat{X}})s^i$
		\EndFor
		\EndIf
		\State \Return $\left(\control + q_1, P_{(\etalb, \etaub)}(\eta + q_2)\right)$
		\EndFunction
	\end{algorithmic}
\end{algorithm}
The use of box-constraints for the extension factor $\eta$ make it necessary to implement the BFGS method similarly to what can be found in \cite{byrd1995limited}, from which \cref{alg::lBFGS} is partly inspired. 
For the box-constrained limited memory BFGS method, we introduce the indicator function for the condition $\etalb \leq \eta \leq \etaub$ as
\begin{equation}
	\chi_\eta(x) := 
	\begin{cases}
		1, \quad \text{ if } \etalb \leq \eta - \sigma \nabla \eta \leq \etaub\\
		0, \quad \text{ else}
	\end{cases}
	\label{eq:eta_indicator}
\end{equation}
for some small $\sigma > 0$.
Recall that the canonical inner product on $X$ is given as 
\begin{equation}
	((\control_1,\eta_1), (\control_2,\eta_2))_X = (\control_1,\control_2)_{L^2(\Gobs)} + (\eta_1, \eta_2)_{L^2(\Omega)}.
\end{equation}
This is now modified to take the active box-constraints into account by introducing $\chi_\eta$ into the second term and thereby reducing the integration to the region of inactive constraints
\begin{equation}
	((\control_1,\eta_1), (\control_2,\eta_2))_{\hat{X}} = (\control_1,\control_2)_{L^2(\Gobs)} + (\eta_1, \chi_\eta \eta_2)_{L^2(\Omega)}
	\label{eq:modified_inner_prod}
\end{equation} 
\Cref{eq:modified_inner_prod} defines the inner product appearing in lines 8, 11, 13 of~\cref{alg::lBFGS}.

Conceptually the optimization scheme presented consists of outer and an inner iterations. 
The outer, seen in~\cref{alg::aug_lag_alg}, updates either $\lambda_g$ or the penalty factor $\tau$ by increment factors $\mult{\mathrm{inc}}$ and $\tau_\mathrm{inc}$, respectively. 
In each cycle of the inner loop, a complete optimization is solved using BFGS updates as seen in~\cref{alg::lBFGS}.
\newline
	\begin{algorithm}[!ht]
		\caption{Computation of Reduced Gradient}
		\label{alg::gradient}
		\begin{algorithmic}[1]
			\Function{Gradient}{$g_\mathrm{def}, \control, \eta$}
			\State $\control \mapsto \defor$ via \cref{eq:weak_deform_equation}
			\State $\defor \mapsto (\vel,p)$ via \cref{eq:ns_variational_1,eq:ns_variational_2}
			\State $(\vel,p,w) \mapsto (\vecmult{\vel}, \mult{p})$ via \cref{eq:weak_adjoint_deform_equation}
			\State $(\vel,p,w,\vecmult{\vel}, \mult{p},\vecmult{g}) \mapsto \vecmult{\defor}$ via \cref{eq:weak_adjoint_deform_equation}
			\State $\vecmult{\defor} \mapsto (\gamma, \kappa)$ via \cref{eq:control_gradient,eq:extension_gradient}
			\State \Return $(w, \gamma, \kappa)$
			\EndFunction
		\end{algorithmic}
	\end{algorithm}
\section{Shape Optimization Applications}\label{sec:results}
In this section we present shape optimization applications with the incompressible, stationary Navier-Stokes equations as state equation.
The purpose of the featured case studies in this section is to show the application of the algorithm presented in \cref{sec:algorithm}, which includes the effect of the nonlinearity control variable $\eta$ on the extension operator $S$.
The obstacle shape deformations demonstrate the algorithm's capabilities at the detection, smoothing and creation of domain singularities such as tips and edges. Aspects of the multigrid preconditioner's effects are discussed.
Moreover, a grid independence study illustrates that the optimal shape is reached regardless of the number of refinements in the grid hierarchy.
The latter result is a fundamental stepping stone towards a scalable parallel implementation of the methodology proposed in \cref{sec:scalability}.

The flow tunnel is depicted as the holdall domain in \cref{fig:2d_holdall} with
\begin{equation*}
	\holdall_{2d} = (-7,7) \times (-3,3) \text{ and } \holdall_{3d} = (-7,7) \times (-3,3) \times (-3,3)
\end{equation*}
for the 2d and 3d cases respectively, taking into account that for the 3d case the obstacle has a spherical shape.
Thus, in 2d we have $\Oobs = (-0.5,0.5)^2$ and $\Oobs = \{ x \in \R^3: \|x\|_2 < 1 \}$ in 3d, respectively.

The boundary conditions at the inflow boundary $\Gin$ are set as 
\begin{equation*}
	\vel_\infty=\left(\max\left\{0, \prod\limits_{i=2}^d \cos(\frac{\pi|x_i|}{\delta})\right\}, 0, \dots, 0\right) \in \R^d
\end{equation*}
with $\delta$ the diameter of the flow tunnel. The side length of the square obstacle is $d = 1$, whereas the radius of the sphere in the 3d case is $r = 0.5$. 
The simulations are performed using UG4 \cite{Vogel2014ug4}. 
We expand UG4 through its C++ based plugin functionality. 
The code used for the studies here presented can be consulted at the online repository in \cite{gmgshapeopt2021}.
The 2d and 3d grids are generated using the GMSH toolbox \cite{GMSH}.
\subsection{2d Results}\label{ss:2d}
In this section we present 2d simulations for a flow with viscosity $\nu=0.03$. 
All PDEs are discretized using a $P_1$ approximation, except for the Navier-Stokes equations and its adjoint which are solved with a stable $P_2-P_1$ finite element discretization. 
For this example, $\eta$ has initial value of 0.5 and box constraints are $0 \leq \eta \leq 1.0$ and $b = 0.001$. 
The grid consists of 421,888 triangular elements, with 5 refinement levels.
\begin{figure}[!htbp]
	\includegraphics[width=0.49\textwidth]{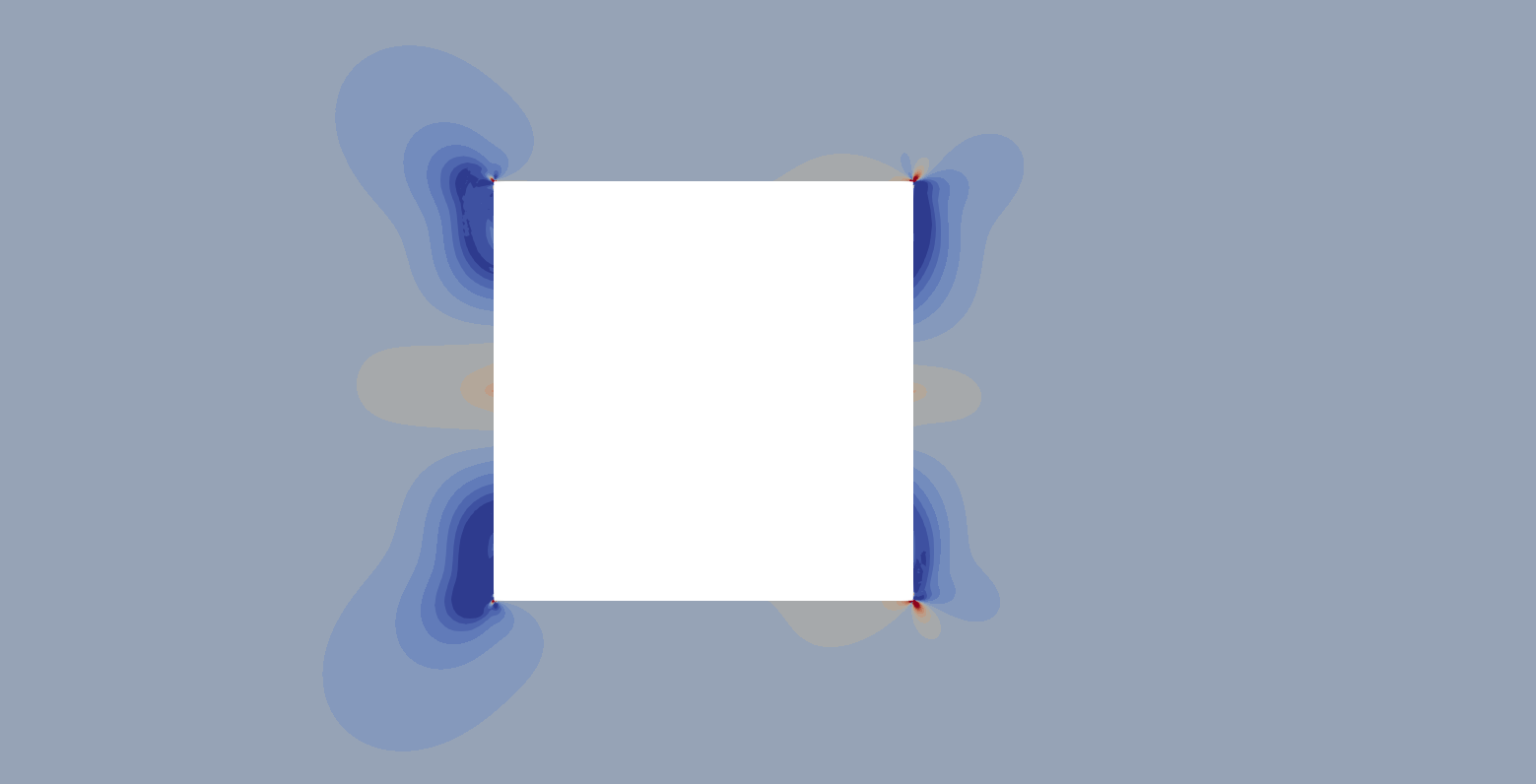} \includegraphics[width=0.49\textwidth]{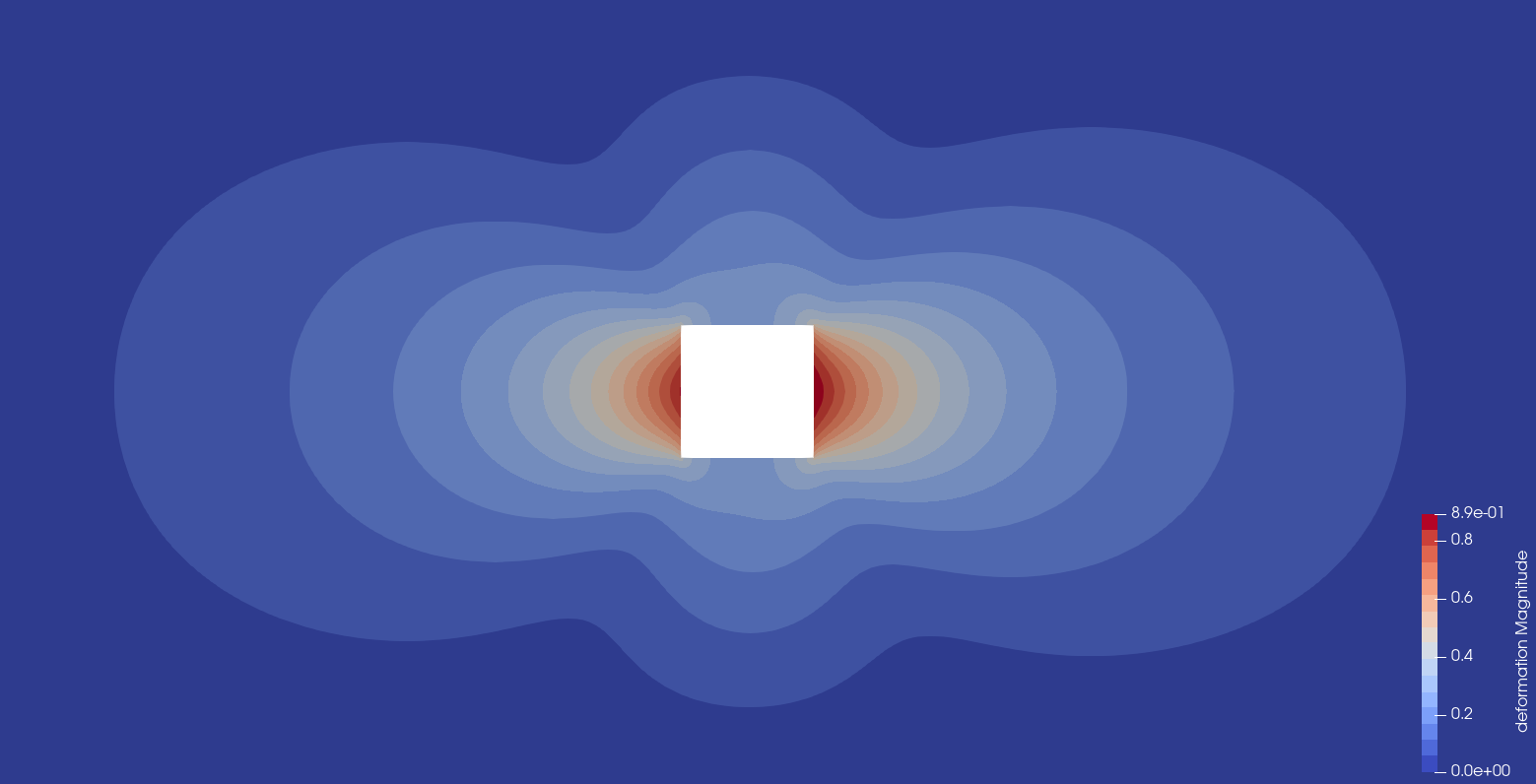} \\\\
	\includegraphics[width=0.49\textwidth]{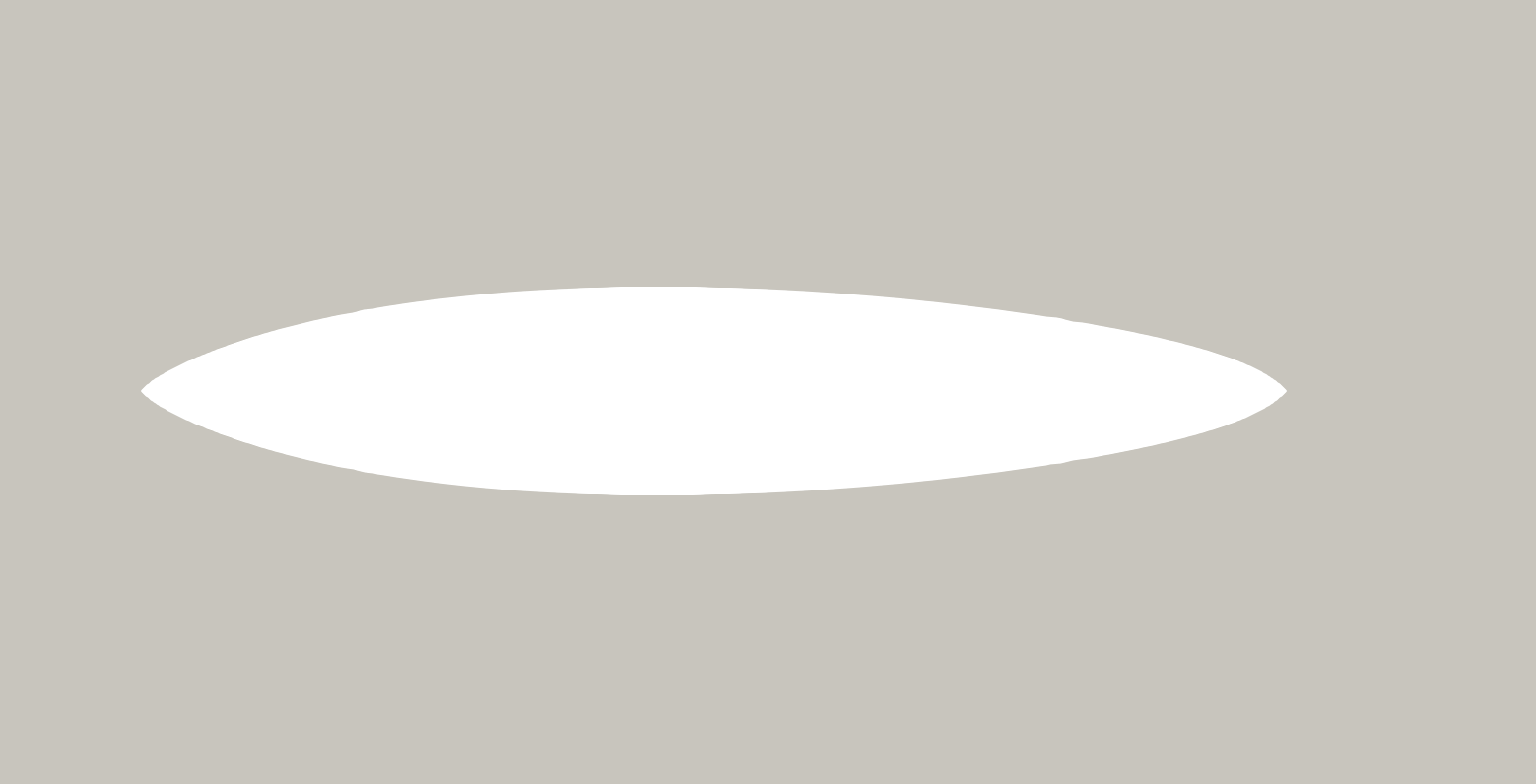} \fbox{\includegraphics[width=0.49\textwidth]{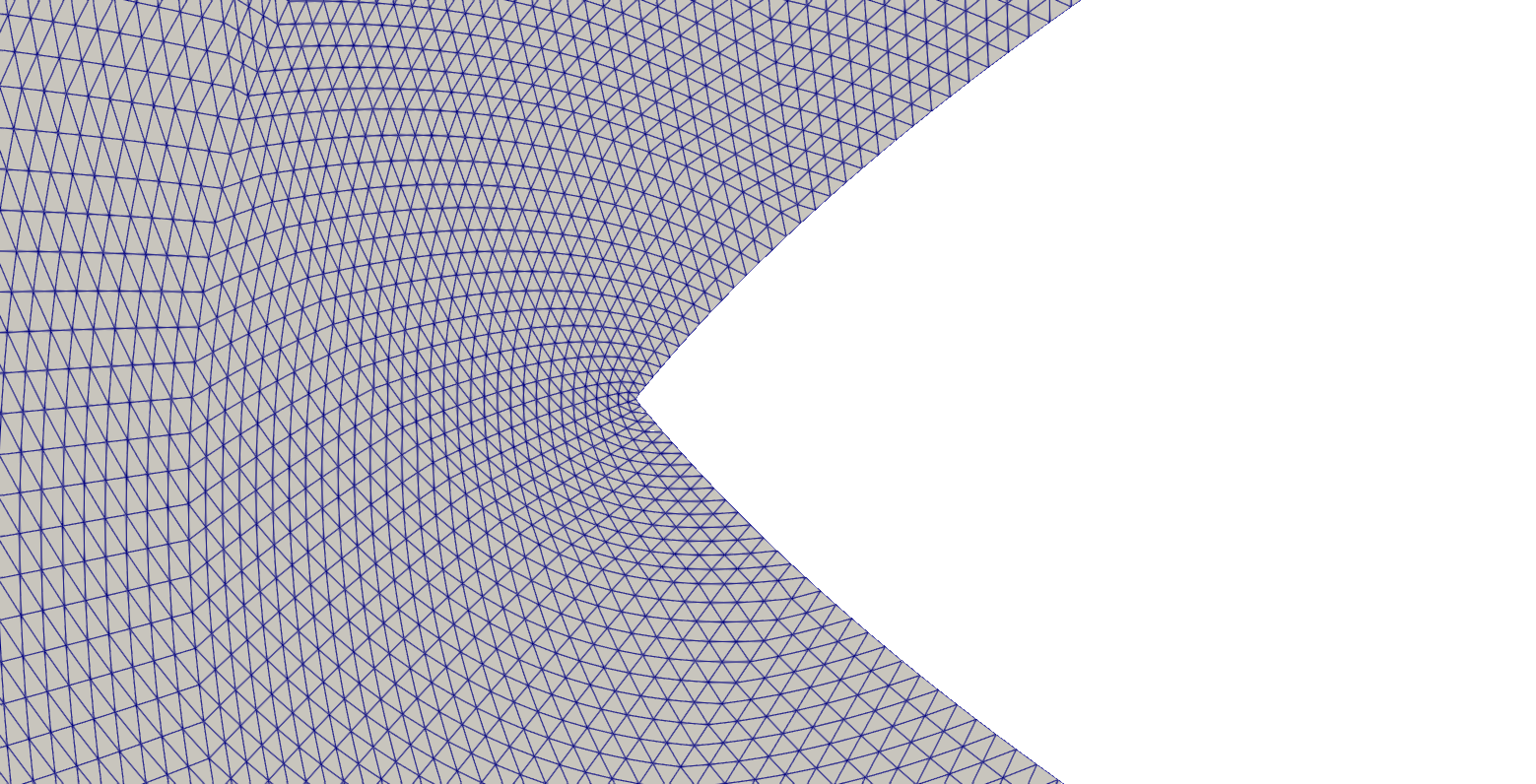}}
	\caption{At the top, the reference configuration is shown with the optimal $\eta$ (left) and $w=S(u)$. At the bottom, the transformed grid $F(\Omega)$ with resulting singularities (left) is shown, altogether with a zoom on the singularity where mesh quality is preserved due to the choice of $S$.}
	\label{fig:2d}
\end{figure}
\Cref{fig:2d} shows results for the optimization of a square obstacle subject to an incompressible, stationary flow. 
The reference configuration with the extension factor $\eta$ and optimal displacement field are shown together with the transformed domain and a closeup of the front tip where the element edges are depicted.
Regarding the reference configuration, it can be seen that the extension factor approaches the imposed values of the box-constraints at two places, the corners of the square and the sections where new singularities have to be created. 
If we recall the weak form of the extension factor \cref{eq:weak_deform_equation}, $\eta$ controls the nonlinearity in each element. 
Given that the same initial value of $\eta$ is set for all elements at the beginning of the simulation and the $\theta$-term in \cref{eq:objective_func} penalizes the deviation from the average of $\tfrac{1}{2}(\etaub+\etalb)$, the different extension factor values, particularly close to the obstacle's surface $\Gobs$, show that equation \cref{eq:weak_deform_equation} adapts depending on the current iterate for the displacement field $\defor$. 
This ensures that the $\defor$ promotes both the generation of new non-smooth points on the boundary, as well as the smoothing of such points introduced by the choice of the reference domain, i.e. the four corners of the box inclusion $\Oobs$. 
This can be observed in \cref{fig:2d}, where high valued displacements are present at the sections where the tips and corners are generated or smoothed. 
Which in turn leads to achieving large deformations $F(\Omega)$ without loss of convergence of the iterative solvers.

In section \cref{sec:background} we already mention that explicit mesh deformations are avoided. 
This comes from the fact that all optimization steps are solved on the reference domain through the method of mappings. 
Therefore we speak of obtaining an optimal deformation field $F$, which is used to transform the domain $\Omega \mapsto F(\Omega)$ according to \cref{eq:opt_sys_2}. 
This is used to obtain the optimal shape, as in \cref{fig:2d}. 
The transformed domain shows the smoothed corners and the generated front and back surface singularities, which are in accordance to the previously mentioned properties of $\eta$ and $\defor$. 
However, throughout the optimization process the proposed algorithm does not require the nodal positions to be redefined, since the reference grid transformation is only performed as part of the post-processing and not of the optimization. 
The close-up corresponds to the front singularity with respect to the direction of flow. 
\Cref{fig:2d} also shows that elements around the generated tip show no distortion and no significant loss of quality. 
This stems from both the effect of the nonlinear term in the extension operator $S$ and the imposed upper bound $b$ on the determinant of the deformation gradient $\Det$, given in \cref{eq:opt_sys_4}. 
The latter condition is what preserves local injectivity, thus avoiding the loss of mesh quality.
In \cref{fig:2d} this is shown as the absence of collapsed or overlapping elements, as is previously mentioned, even for the elements that clearly undergo large deformations, i.e. the ones that conform the generated tips and the smoothed square corners. 
Moreover, in \cref{sec:scalability,ss:grid_stdy} this can be understood as a mesh independent preservation of the geometrical and numerical convergence in terms of the final optimal shape achieved and the total iteration counts of the iterative solvers.

\begin{figure}[!htbp]
	\centering
	\begin{tabular}{rccccc}
		Step & \hspace{0.1em} & Transformed domain & Extension factor \\
		${2}~$ &  
		&\includegraphics[width=0.45\textwidth]{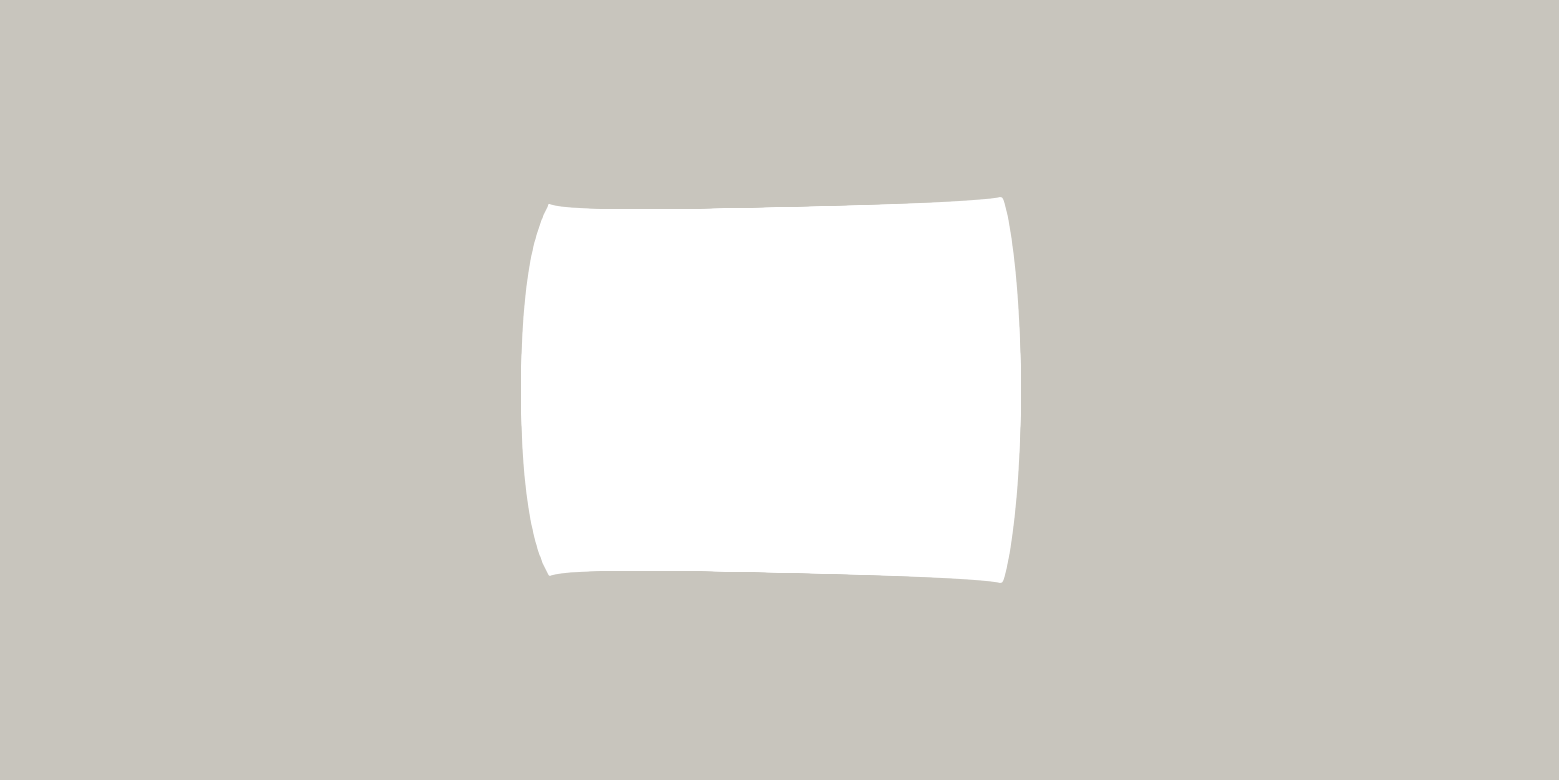} &
		\includegraphics[width=0.45\textwidth]{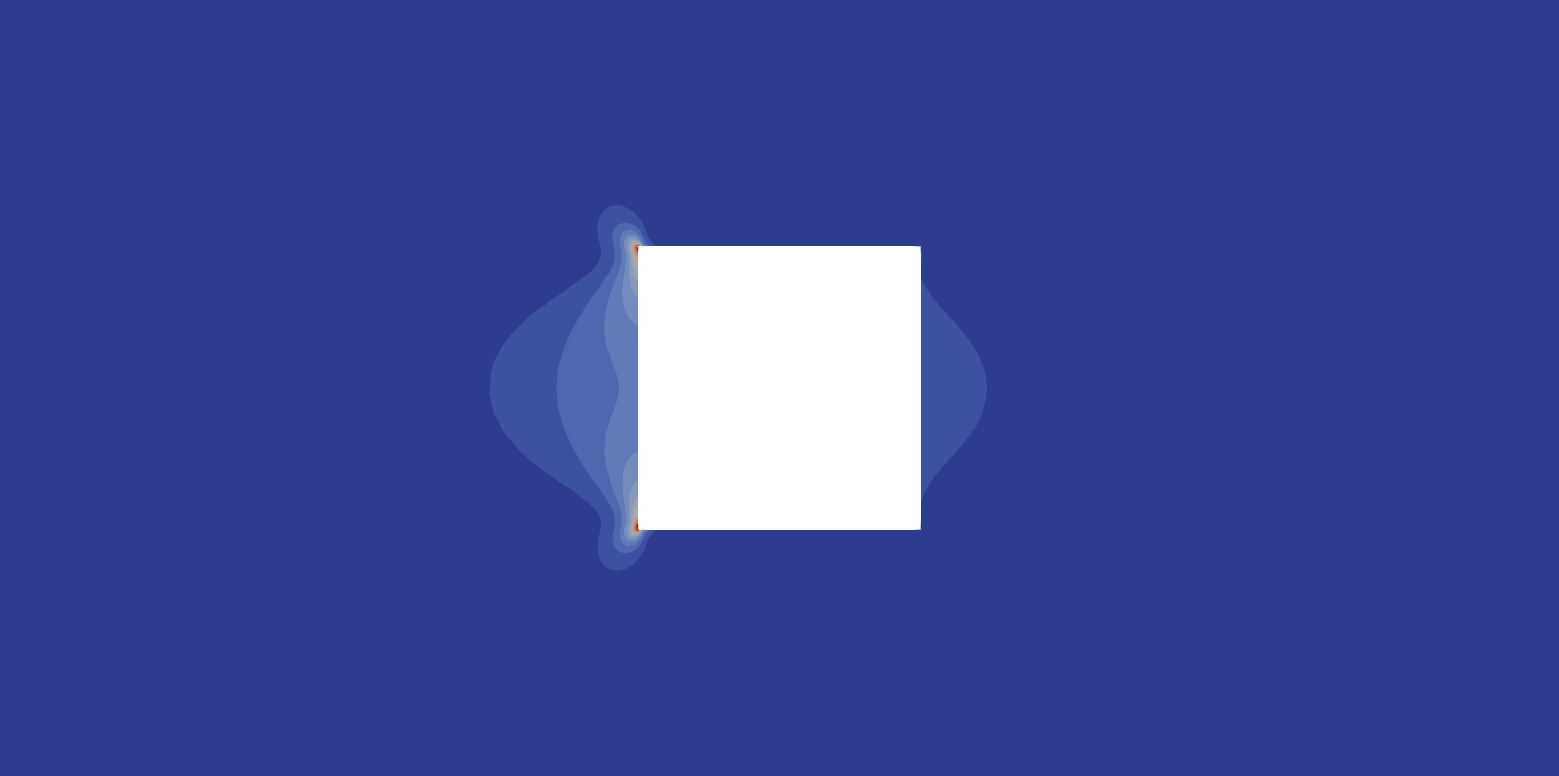} \\
		${8}~$&
		&\includegraphics[width=0.45\textwidth]{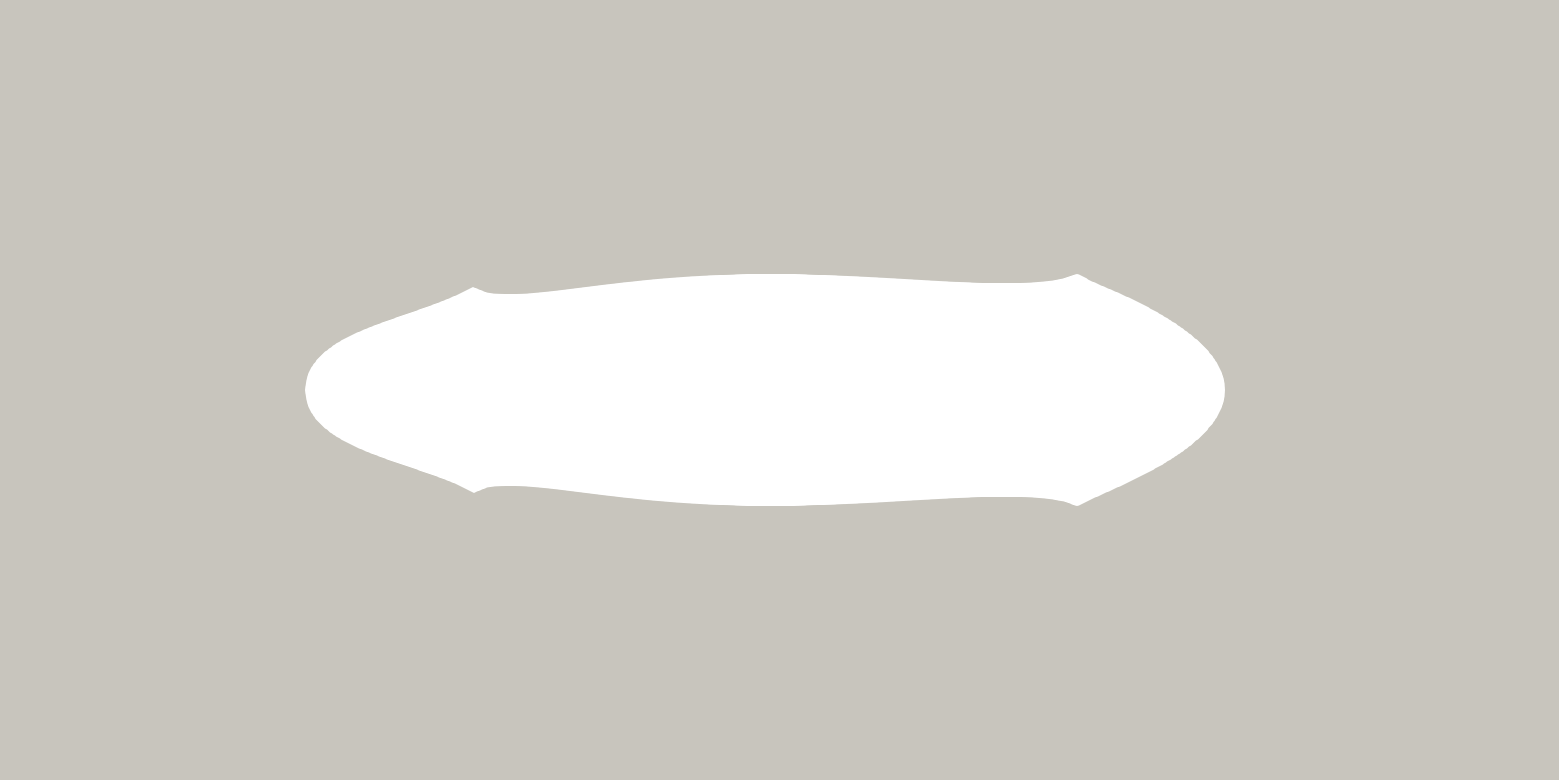} &
		\includegraphics[width=0.45\textwidth]{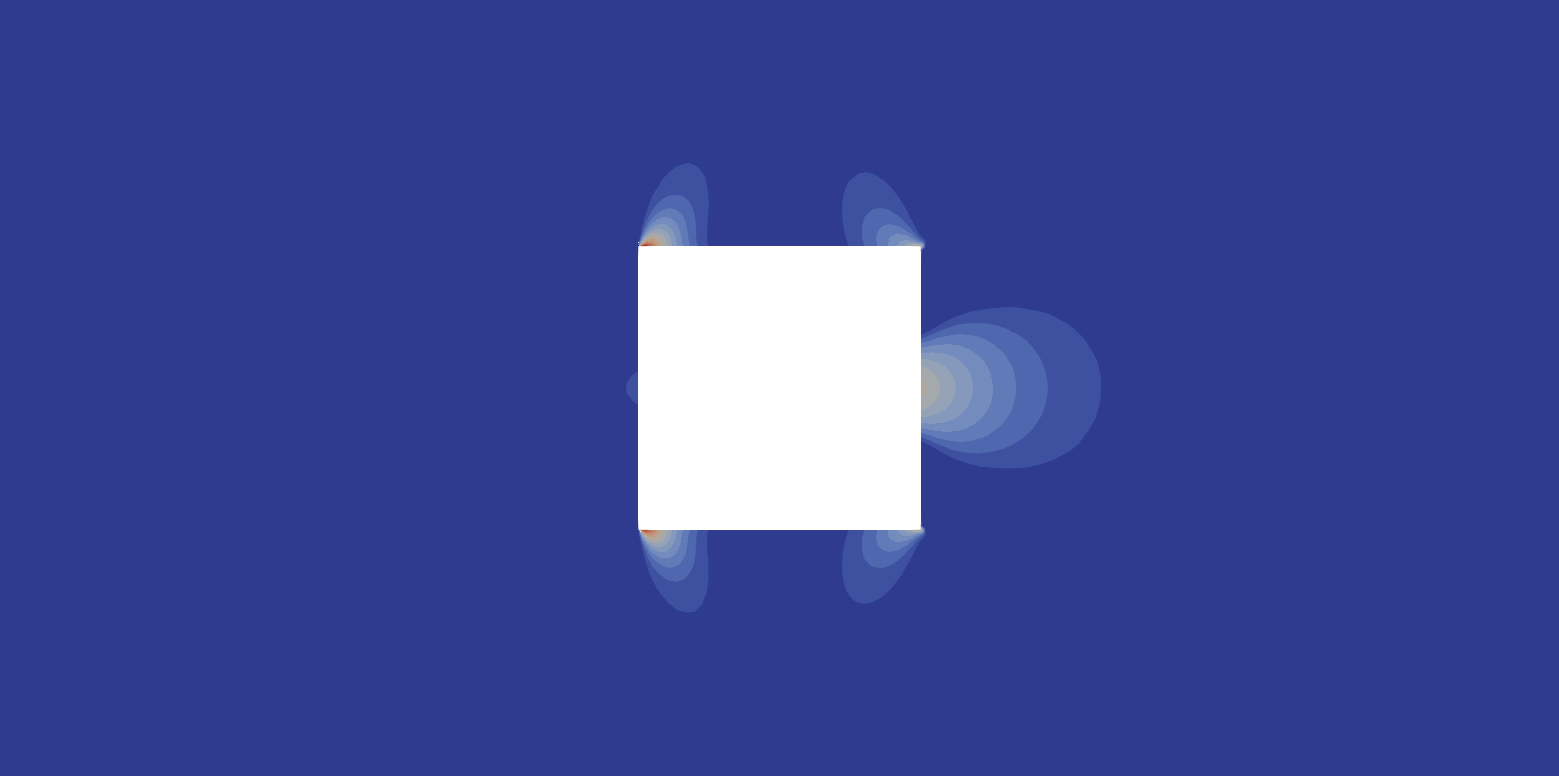} \\
		${20}~$ &  
		&\includegraphics[width=0.45\textwidth]{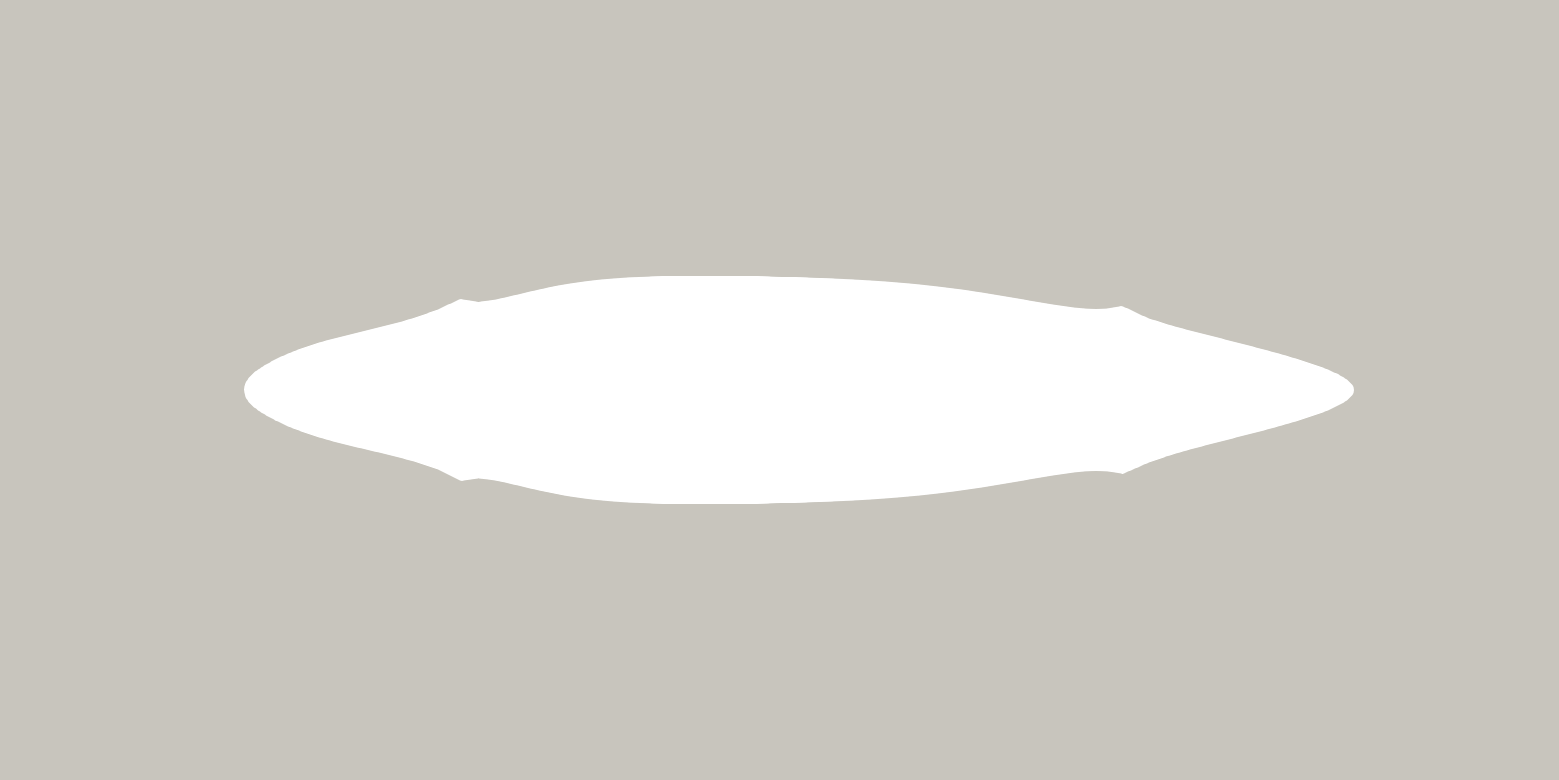} &
		\includegraphics[width=0.45\textwidth]{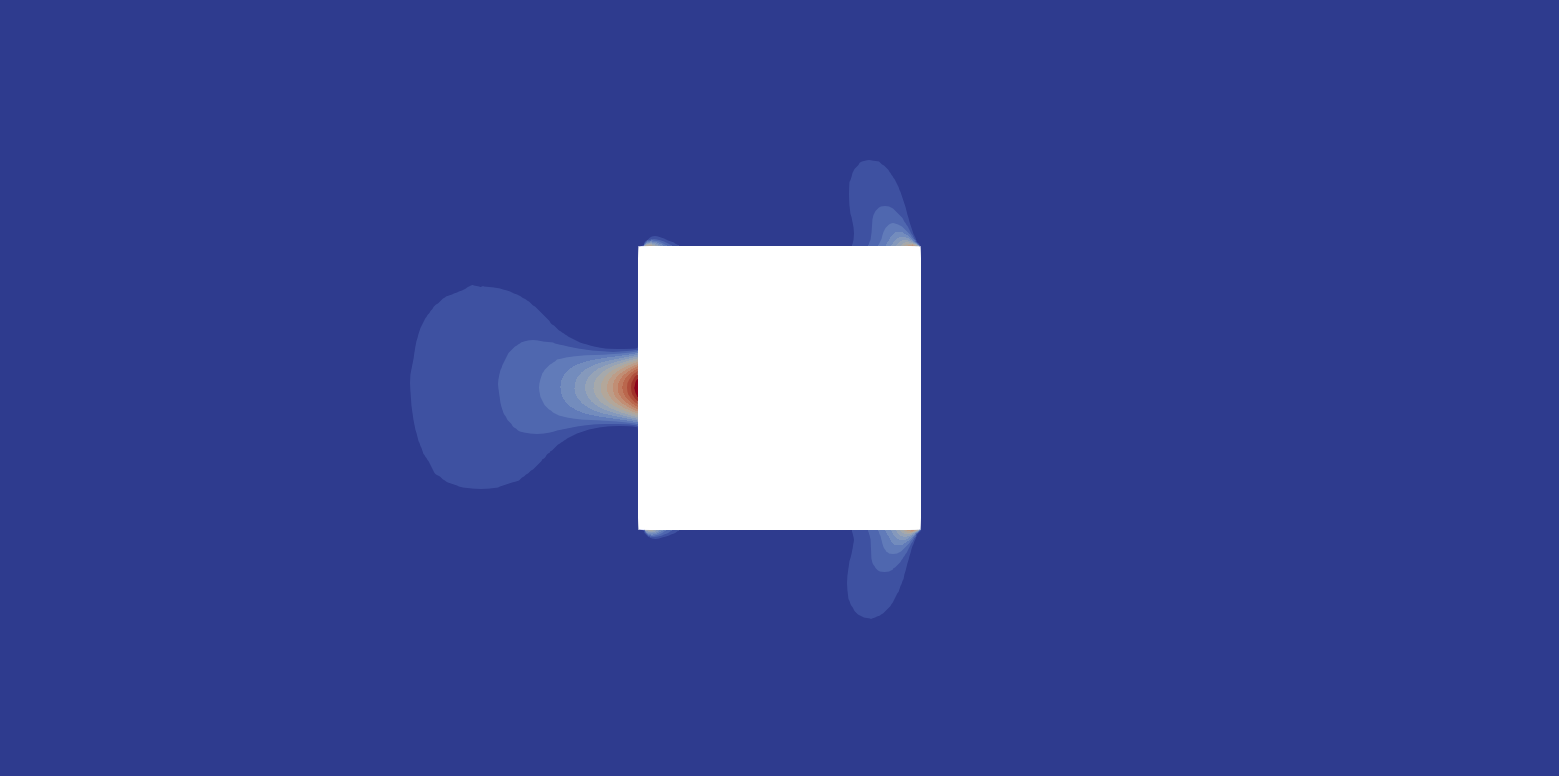} \\
		${56}~$ &  
		&\includegraphics[width=0.45\textwidth]{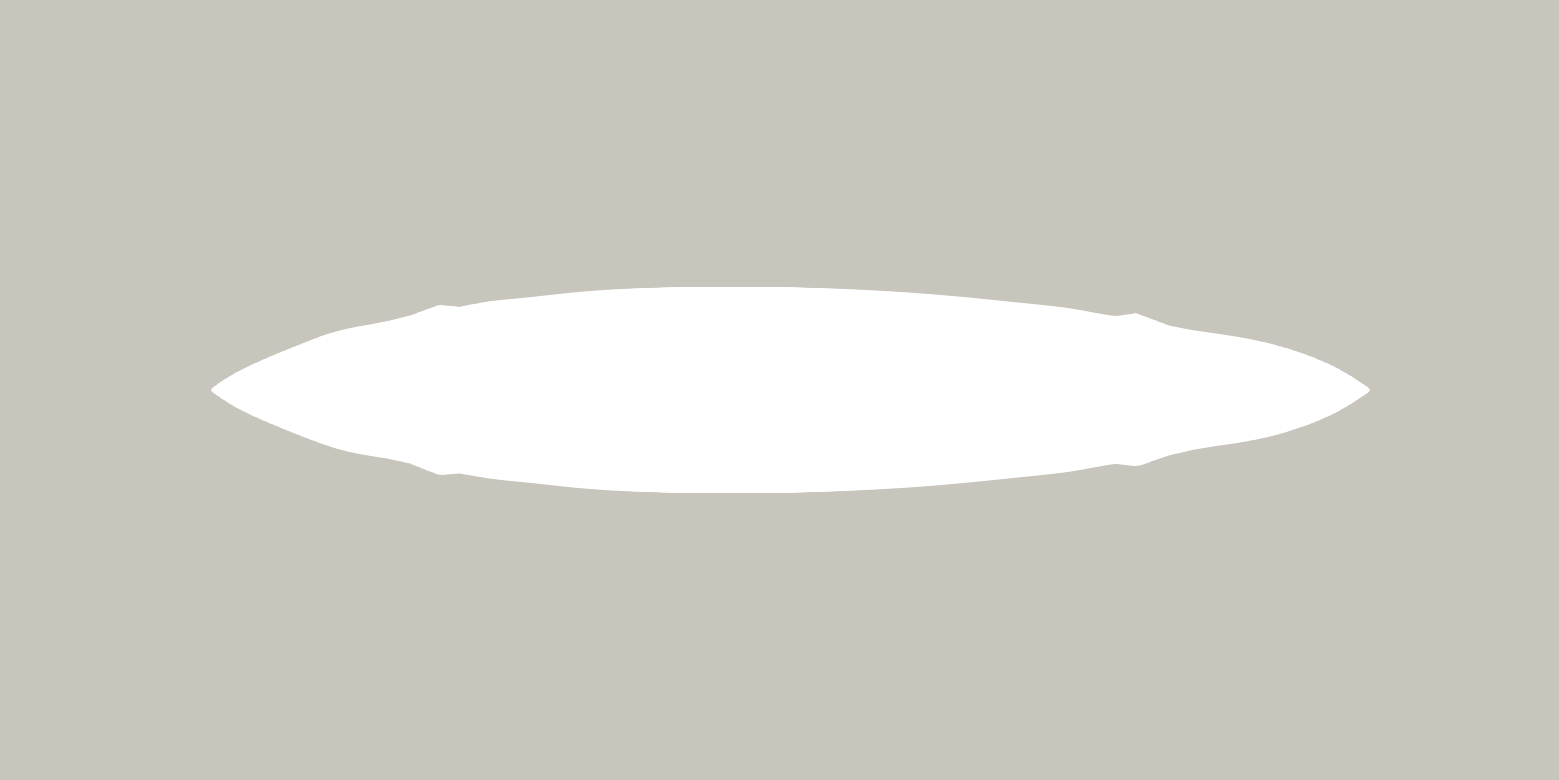} &
		\includegraphics[width=0.45\textwidth]{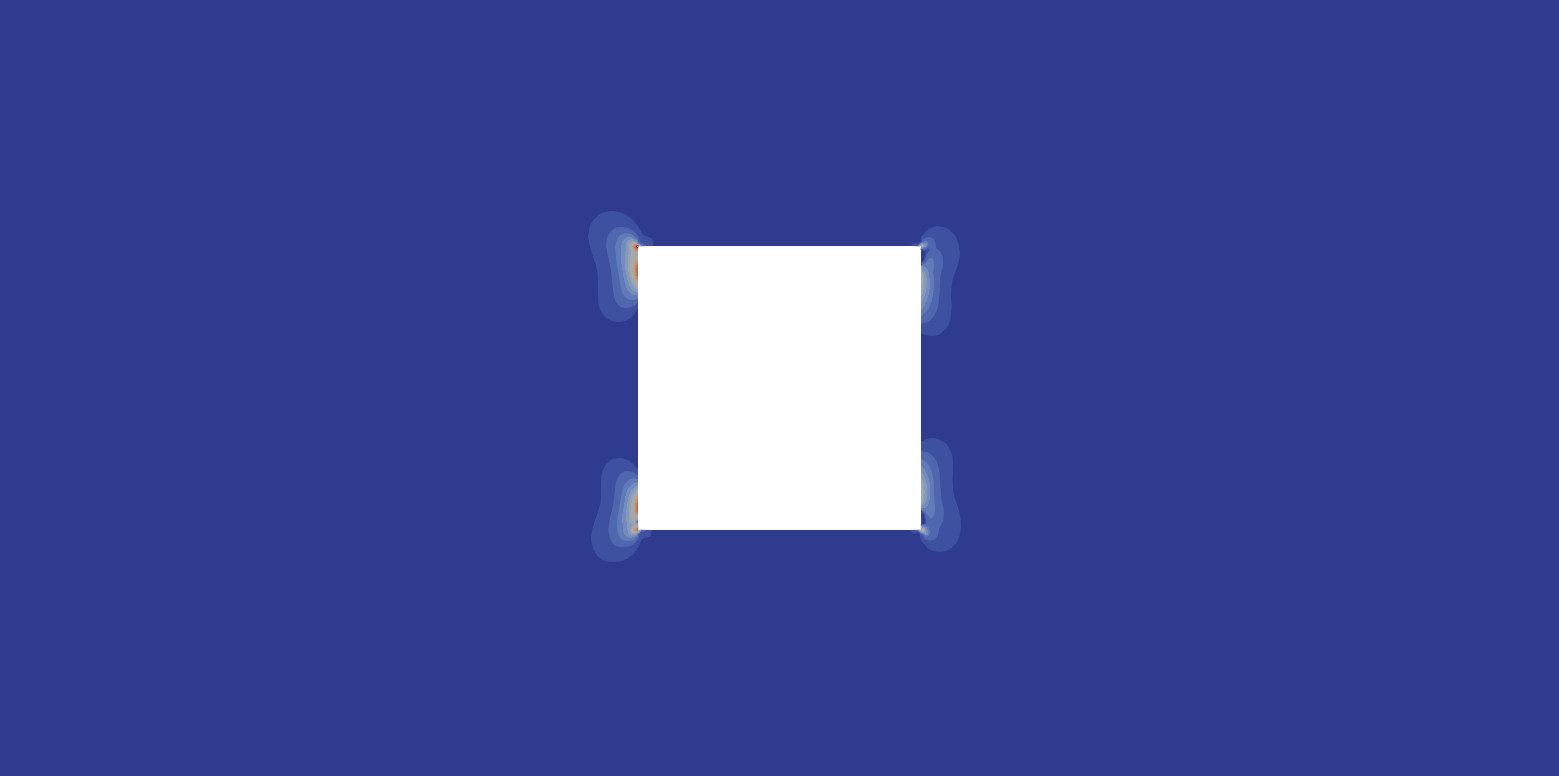} \\
		${74}~$ &  
		&\includegraphics[width=0.45\textwidth]{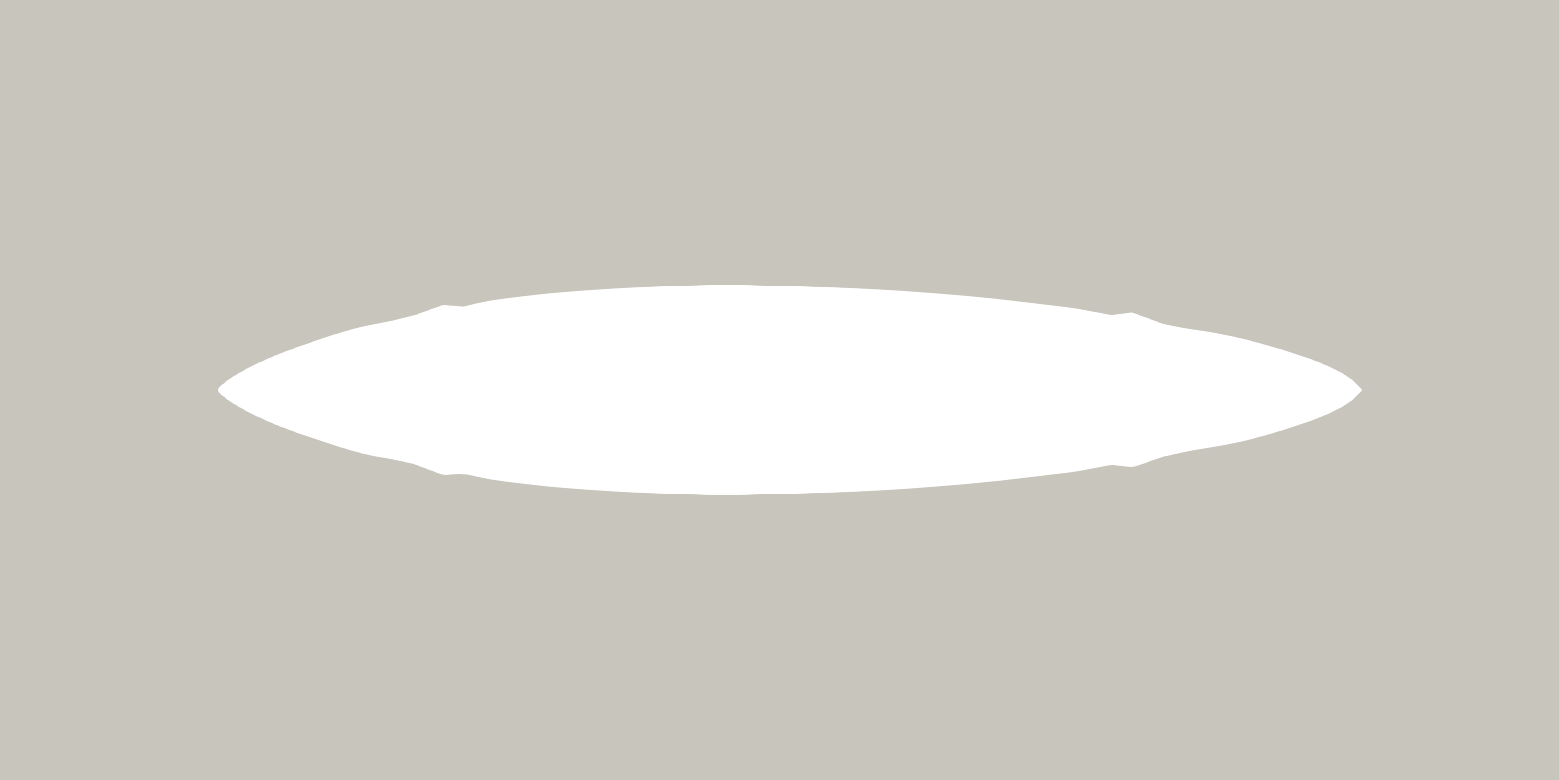} &
		\includegraphics[width=0.45\textwidth]{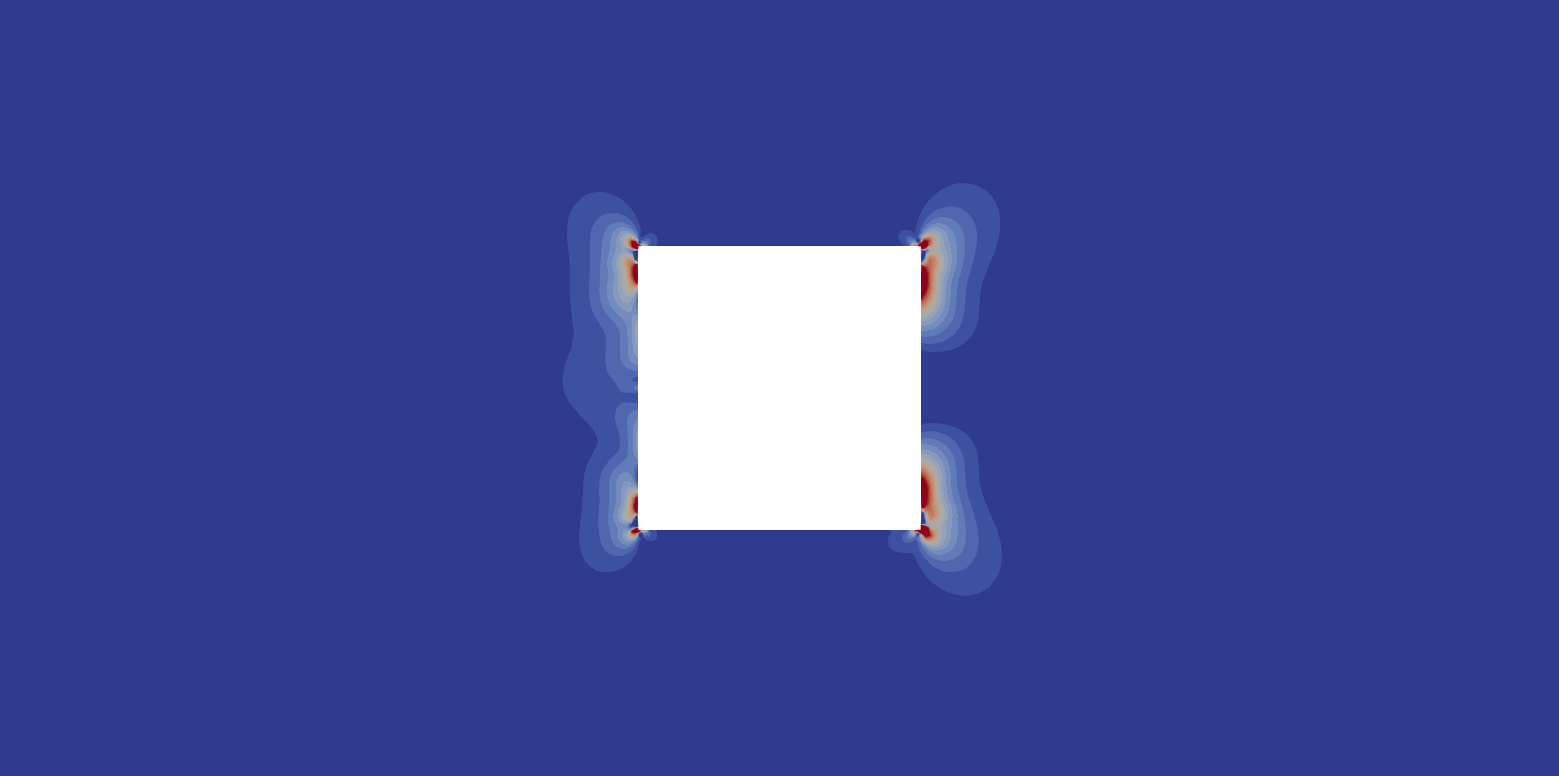} \\
	\end{tabular}
	\caption{Transformation of the reference domain by the application of the deformation field compared to the accumulation of the extension factor, given for steps 2,8,20,56,74}
	\label{fig:2dDeformationSequence}
\end{figure}

On the other hand, the extension equation adapts where the tips have to be created to reach an optimal value of the objective function.
This is illustrated through the changing value of the nonlinearity switch $\eta$ on each optimization step. 
\Cref{fig:2dDeformationSequence} shows the plot of $\eta$ over the reference domain compared to the domain transformed by the displacement field $\defor$.

At the start of the simulation the extension equation is already adapted to find the corners of the reference configuration, this is shown as the concentrated value of $\eta$. 
As the obstacle's initial singularities are removed, necessary ones are created. 
This causes a concentration of the extension factor at the reference configuration locations where new geometrical singularities have to be created. 
Afterwards, the optimization scheme works towards smoothing the obstacle's surface, therefore $\eta$ goes through no major concentration values across the grid, as can be seen in step 74 of~\cref{fig:2dDeformationSequence}.

The distribution of $\eta$ across the grid has to be compared against the transformed grid.  
Given that the Lagrange multipliers are yet to converge, the initial steps can incur in violations of the geometrical constraints. This can be seen as the highly deformed shapes at the initial steps of~\cref{fig:2dDeformationSequence}. However, as the algorithm performs the multipliers' update, as in~\cref{alg::aug_lag_alg} line 12, the geometrical constraints are fulfilled according to the prescribed $\epsilon_g$ and the new obstacle surface's singularities are formed. Moreover, since the reference configuration singularities are identified at the initial optimization steps, the necessary smoothing is carried out until the simulation converges or the maximum number of steps is reached. This can be seen comparing the last 2 steps of~\cref{fig:2dDeformationSequence}.
\subsection{3d Results}\label{ss:3d}
\begin{figure}[!htbp]
	\centering
		\begin{minipage}{0.5\textwidth}
			\includegraphics[width=1.0\textwidth]{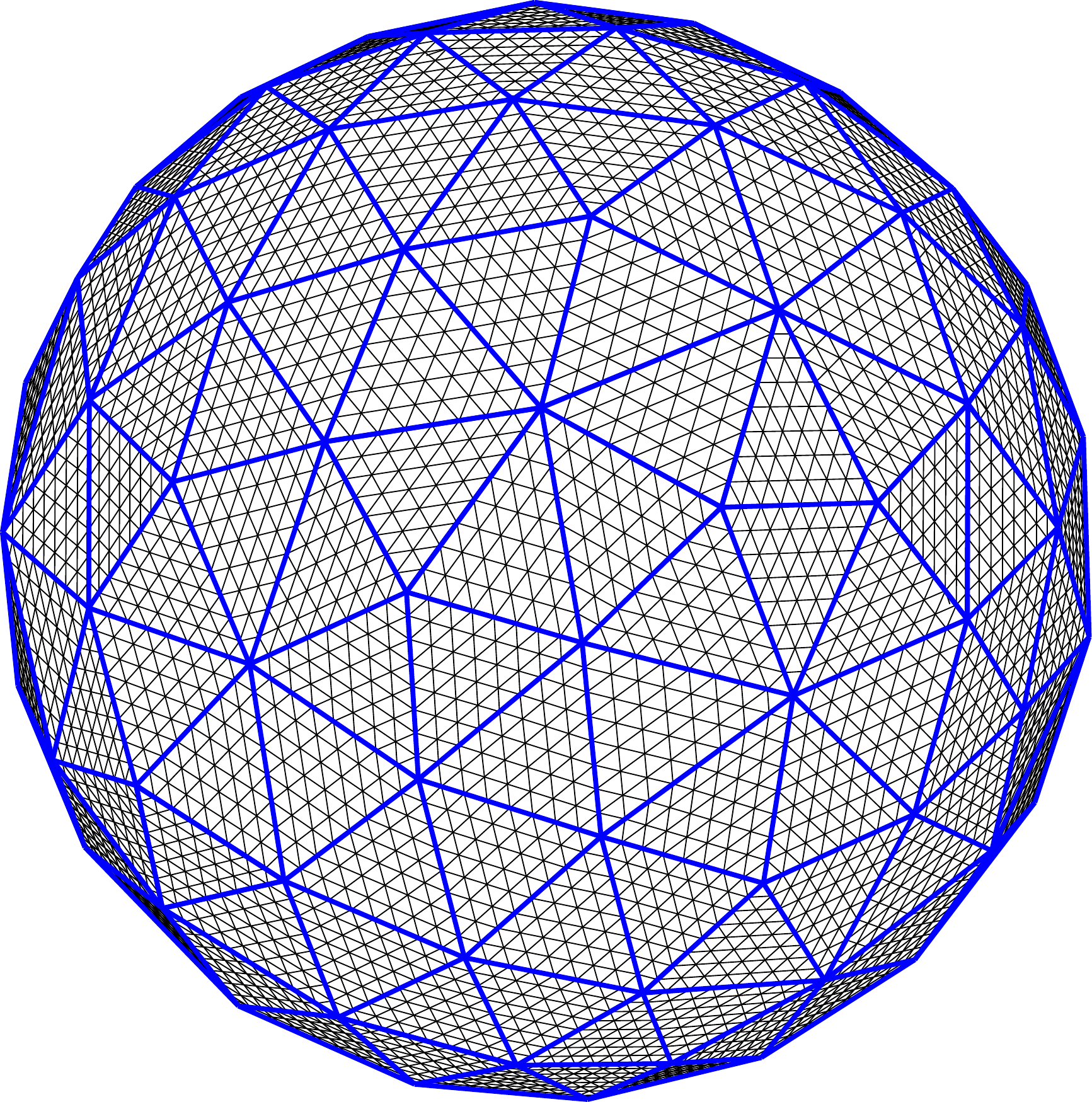}
		\end{minipage}
	\caption{3d highest grid level is shown compared to the base level (bold lines) for the unit-diameter sphere obstacle.}
	\label{fig:3d_grid}
\end{figure}
3d results for the optimization of a unit-diameter sphere are presented here. 
For these results, 4 levels of refinements are used with up to 12,296,192 tetrahedral elements, while the obstacle's surface consists of 54,784 triangular elements. 
The viscosity is set to $\nu=0.1$, with a discretization scheme as in \cref{ss:2d} where, in contrast to the 2d case, $P_1-P_1$ mixed elements are used for the Navier-Stokes equations and its adjoint.
Regarding the extension equation, $\eta$ has initial value of 30 and the box constraints are $0 \leq \eta \leq 60.0$.
A pressure projection stabilization term is used for the mixed finite element approximation, as given in \cite{elman2014}.

The discretized domain representing the 3d obstacle is shown in \cref{fig:3d_grid}. 
Both the base level (bold blue lines) and highest refinement level are presented. 
As mentioned in \cref{sec:intro}, we investigate the application of the geometric multigrid method as a preconditioner in shape optimization. 
This implies that we strive to maintain the base level as coarse as possible, as can be seen by the underrepresented sphere shown; with the idea of solving the coarsest problem with a direct method as quickly as possible. 
While this is ideal for the usage and convergence of the geometric multigrid method, it has some undesired effects. 
As can be seen, the several refinements introduced by the creation of the hierarchical grid levels do not correct necessarily introduce a smoothing of the obstacle's surface. 
The refinements are limited to subdividing the triangular faces present on $\Gobs$, while the edges from the base grid remain.
\begin{figure}[!htbp]
	\centering
	\begin{minipage}{0.49\textwidth}
		\includegraphics[width=1.0\textwidth]{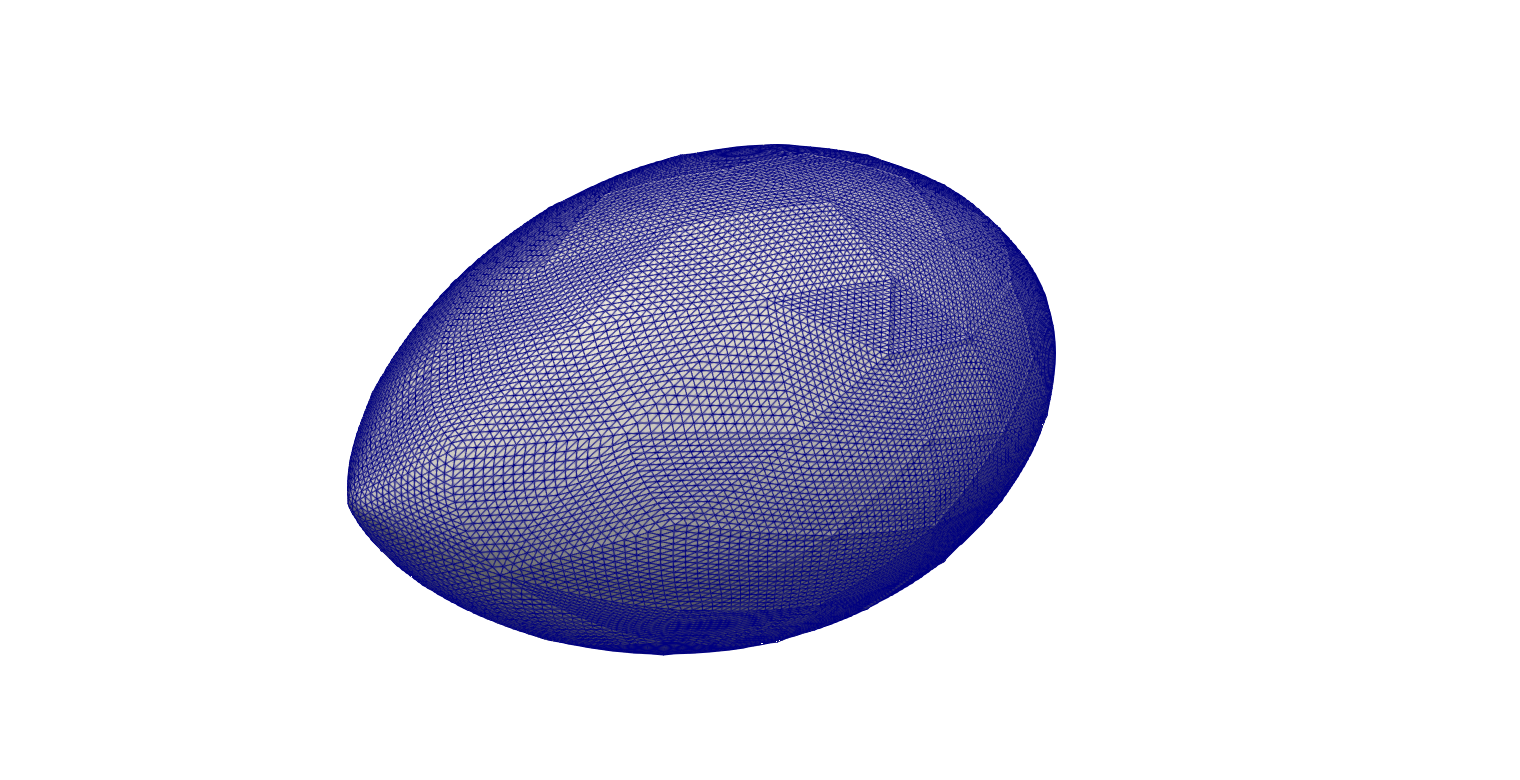}\\
		\includegraphics[width=1.0\textwidth]{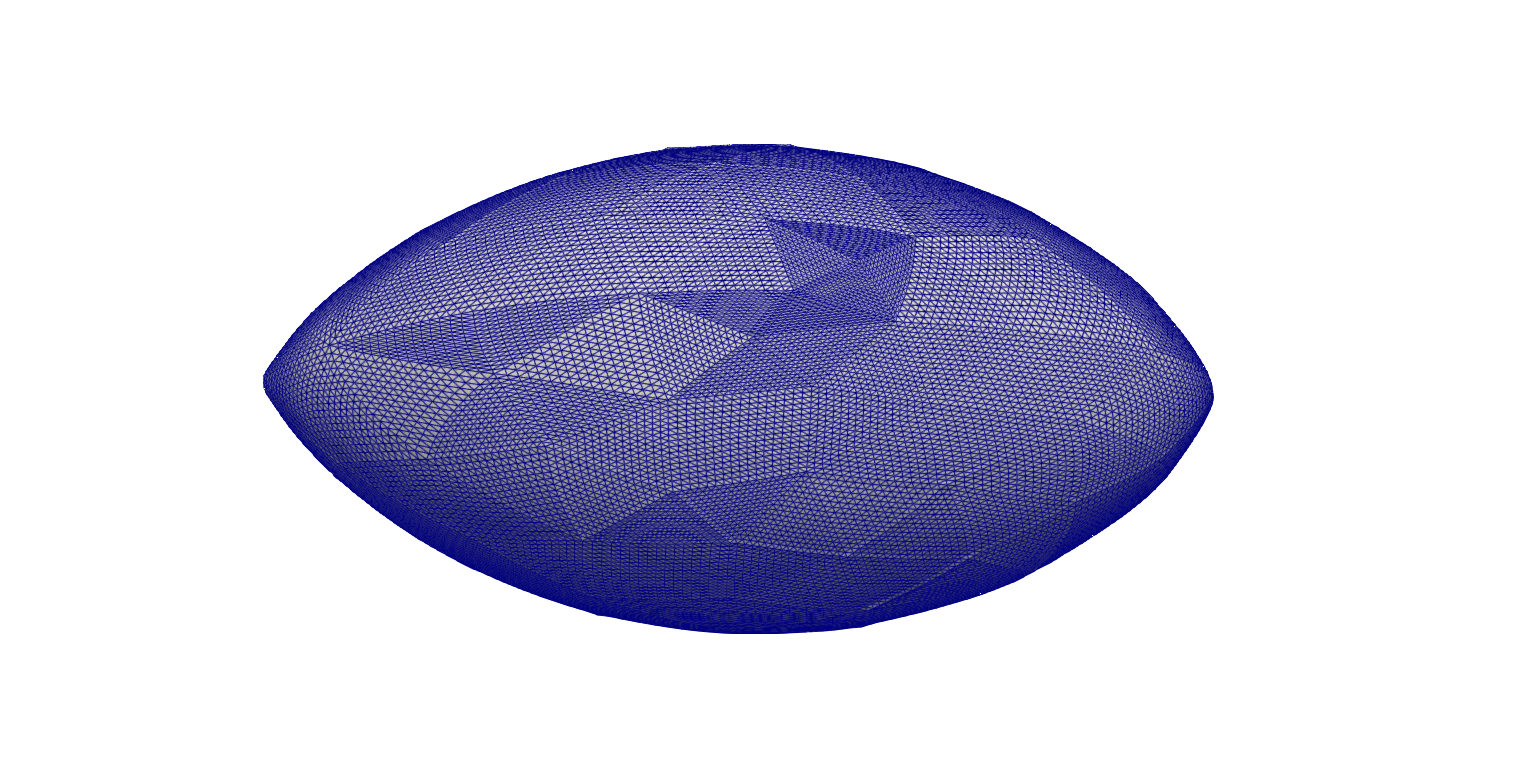}
	\end{minipage}
	\caption{3d optimal shape on the highest level of refinement in flow with large viscosity}
	\label{fig:3d}
\end{figure}

The results after 61 steps are shown in \cref{fig:3d}.
Non-smooth points, i.e.\ the two tips, are generated on the front and back of $\Gobs$ with respect to the direction of flow.
This is comparable to the optimal shape obtained for the 2d case in \cref{fig:2d}.
The effects of the grid hierarchy can be seen as the remaining edges of the super-elements.
\subsection{Grid Independence Study}\label{ss:grid_stdy}
In order for the proposed optimization scheme to be scalable in terms of time-to-solution to very high numbers of DoFs, it is necessary for the obtained obstacle shape to be independent of the initial level of refinement.
In other words, besides the scalability of the finite element building blocks of the optimization algorithm, the overall convergence of the objective function has to be mesh independent.
This can be understood as obtaining the same optimized shape after a given number of outer iterations in \cref{alg::aug_lag_alg}, with the necessary surface singularities appearing at approximately the same locations. 
Therefore, in this section we provide results for a comparative study between different levels of refinement.
\begin{figure}[!htbp]
	\centering
	\includegraphics[width=0.8\textwidth]{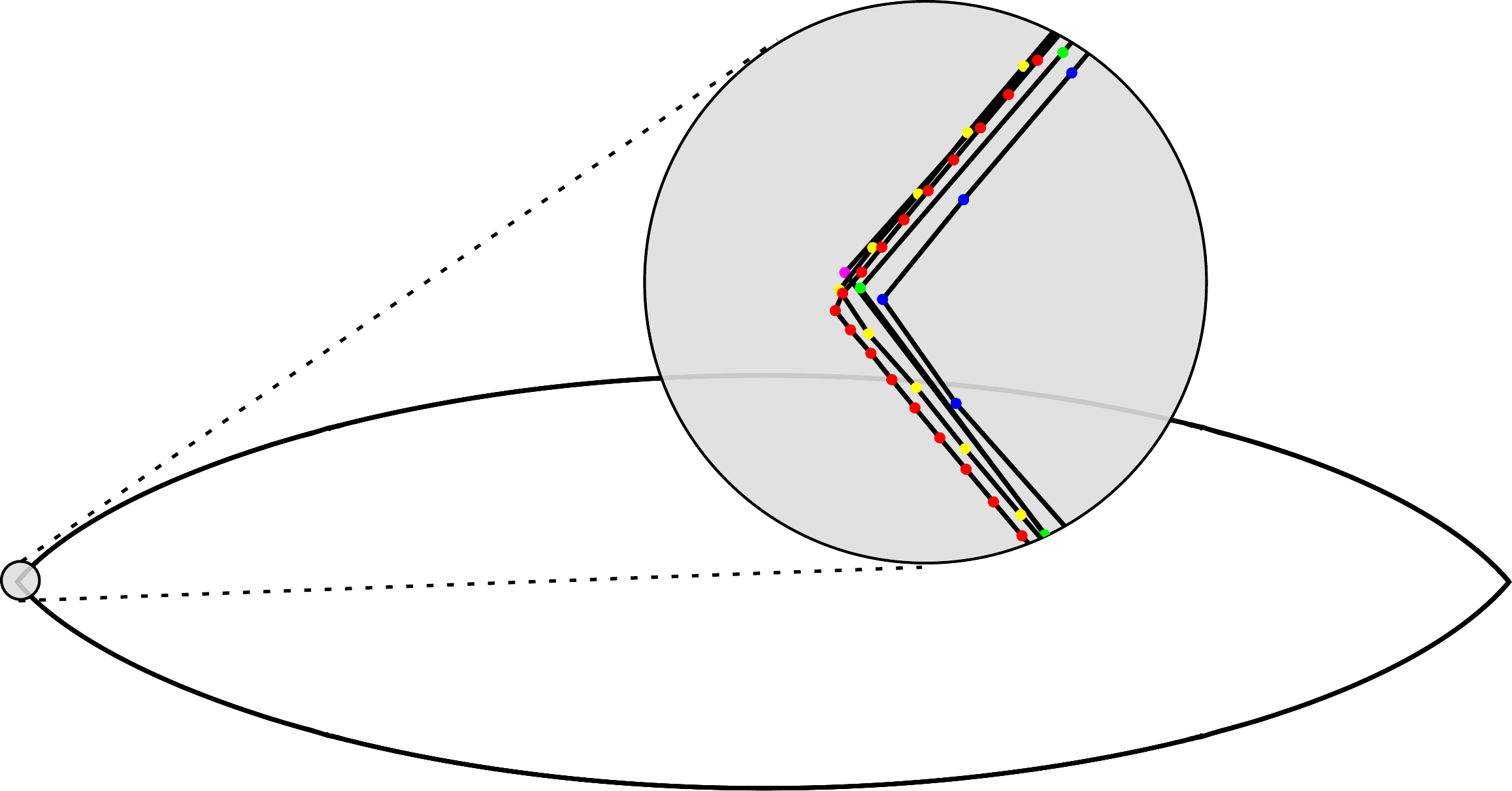}
	\caption{Optimal displacement field $\defor$ after 400 optimization steps applied to the reference shape for several levels of refinement, indicated by colored nodes.}
	\label{fig:gs_shapes}
\end{figure}

\begin{figure}[!htbp]
	\centering
	\includegraphics[width=0.6\textwidth]{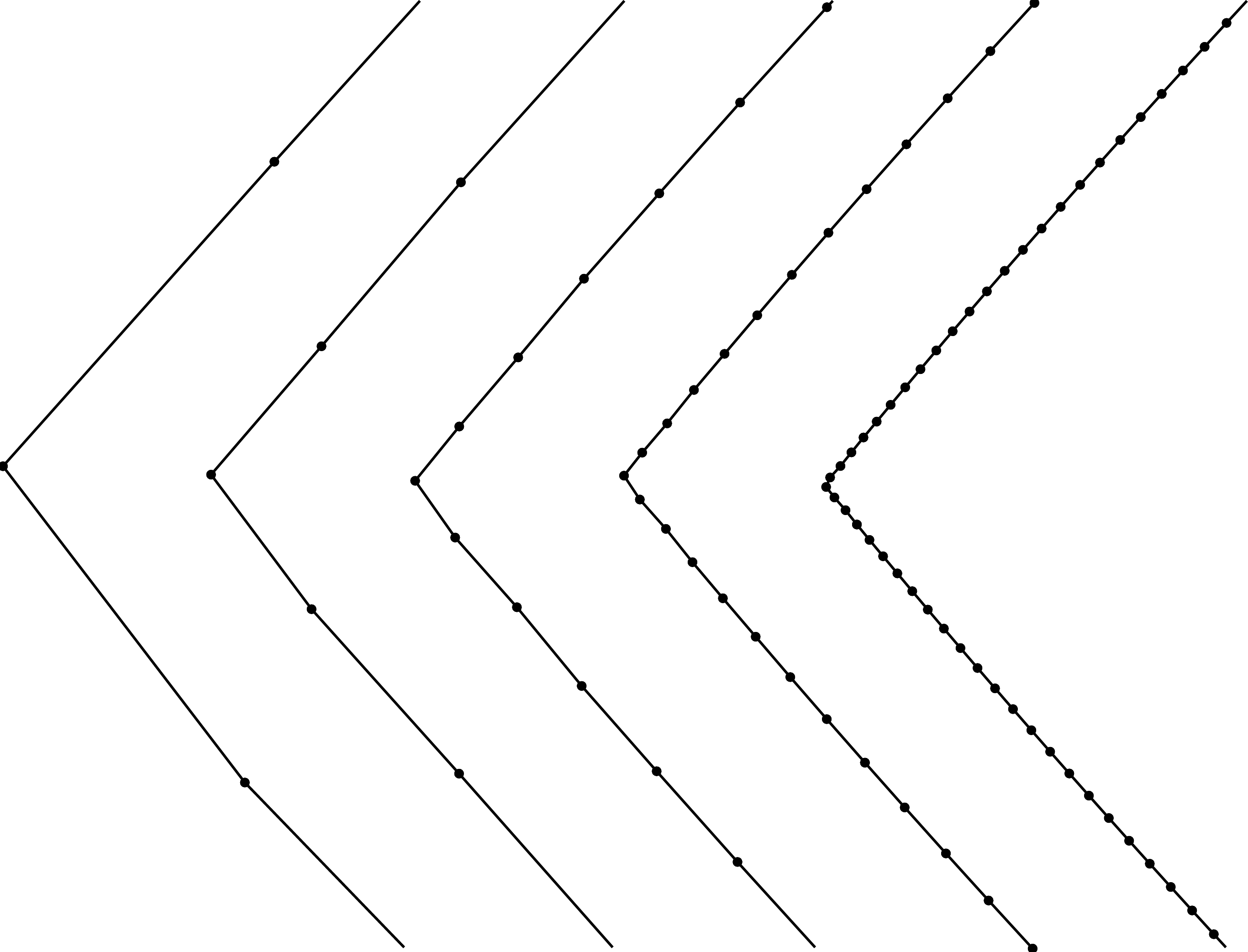}
	\caption{Zoom into the tips of the deformed shapes of the experiment shown in \cref{fig:gs_shapes}. Comparison of generated boundary singularities for 2,3,4,5,6 levels of refinement (from left to right). Shapes are aligned with fixed interspace for better comparability.}
	\label{fig:gs_tips}
\end{figure}

The grid used in this section is formed by 412 triangular elements and refined to 1,687,552 elements.  
The results shown go from 2 to 6 refinement levels. 
The simulations are set with viscosity $\nu=0.1$. 
An equal number of 400 optimization steps is run for all grids, to have an adequate comparison point.
\Cref{fig:gs_shapes} shows the superimposed contours of the obstacle for 2,3,4,5,6 levels of refinement, indicated by the colored nodes. 
A magnification is used on the front tip, to emphasize that all tips appear on the same location, with slight differences owing to the discretization error introduced by different element sizes. 
In addition to this results, \cref{fig:gs_tips} shows a side-by-side comparison of the tips of the aforementioned refinement levels. 
This indicates, that the singularities on the obstacle surface generated by algorithm \cref{alg::aug_lag_alg}, are grid independent.\\
\begin{figure}[!htbp]
	\begin{subfigure}[t]{0.99\textwidth}
		\centering
		\vspace*{0.1cm}
		\scalebox{0.8}{\includegraphics{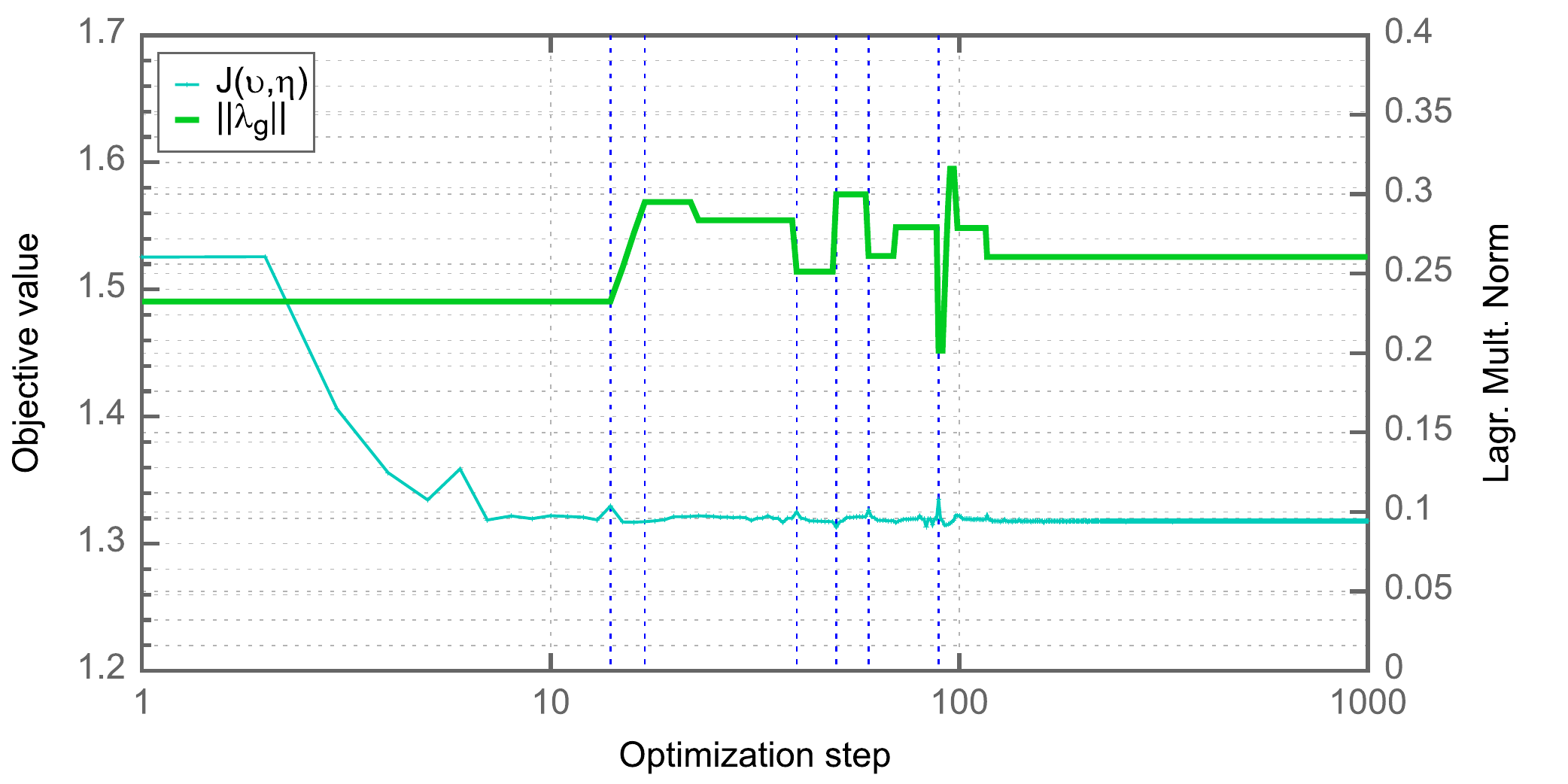}}
	\end{subfigure}
	\begin{subfigure}[t]{0.9\textwidth}
		\centering
		\vspace*{0.1cm}
		\scalebox{0.8}{\includegraphics{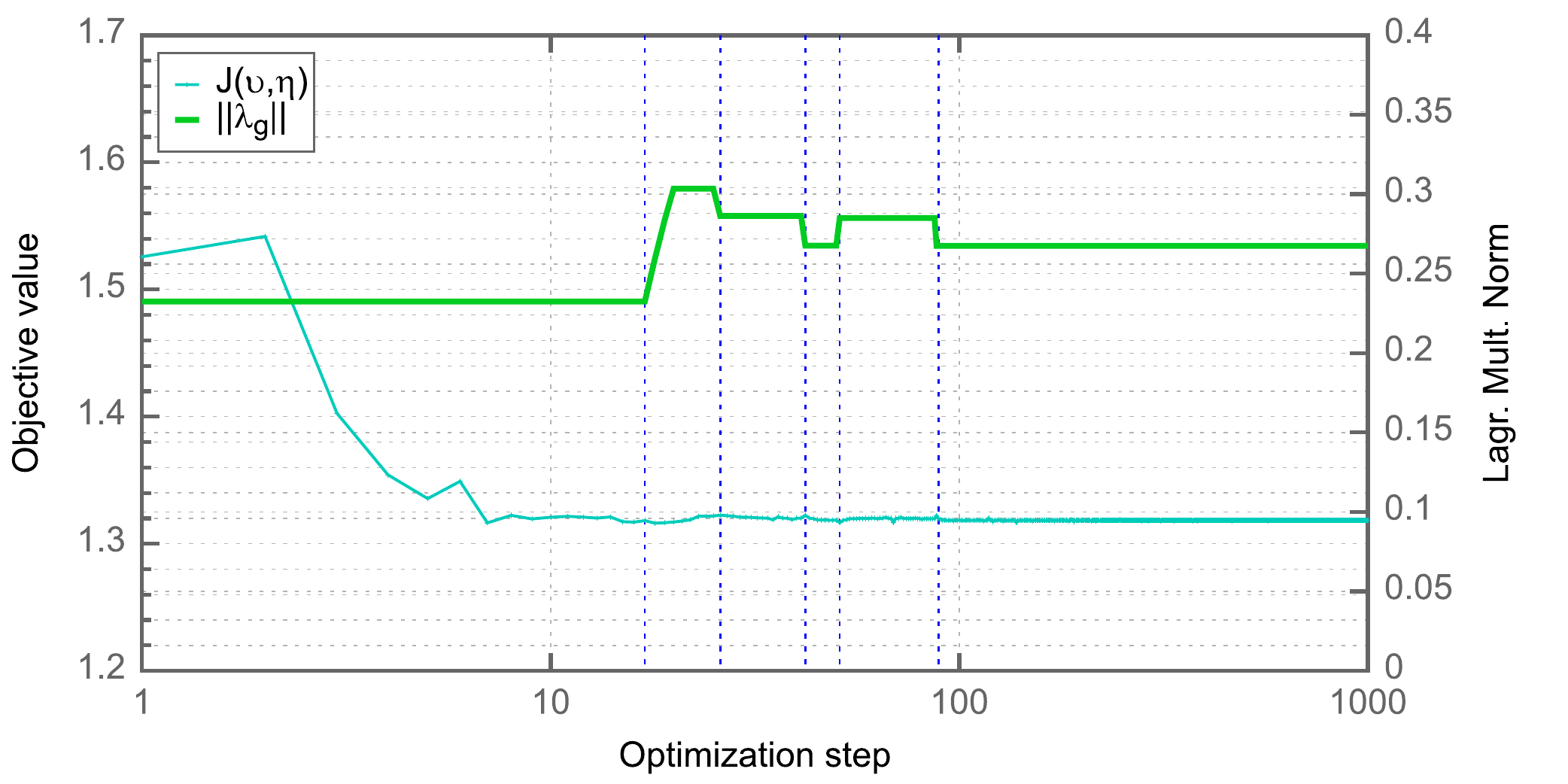}}
	\end{subfigure}
	\caption{Objective function plot for 3 (top) and 4 (bottom) refinements. Green line shows norm of geometric Lagrange multipliers. Dashed blue vertical lines indicate $\lambda_g$ update in augmented Lagrange algorithm.}
	\label{fig:objective_function}
\end{figure}
Moreover, in \cref{fig:objective_function} the objective function plots for 2 different refinement levels are shown. 
With the same viscosity as in \cref{fig:gs_tips,fig:gs_shapes}, the simulations are set to run for 1000 steps. 
The purpose is to demonstrate that the achieved minimal value is independent of the geometry.
We thus choose this large number of optimization steps regardless of any tolerance $\epsilon_\mathrm{outer}$.
The value of \cref{eq:augmented_objective_func} (in blue)  is compared against the Euclidean norm of the Lagrange multipliers (in green) for each refinement, while the update of the multipliers is signaled by the dashed lines (dark blue). 
It is evident that before the convergence of the multipliers, the optimization process is local, which is why differences between the two plots with respect to the objective function value are present.
In \cref{alg::aug_lag_alg}, the condition for the update of the aforementioned Lagrange multipliers is mentioned.
This is related to the set tolerance $\epsilon_g$.
Given that the two geometries are different due to the level of refinement, the fulfillment of the geometrical tolerance is not necessarily achieved in the same optimization steps.
Which in turn, as seen in \cref{alg::aug_lag_alg}, has an effect as to when the multipliers are updated.
Nevertheless, as previously mentioned in this section, the objective function converges altogether with the norm of the multipliers, as seen by comparing the plots in both refinement cases presented in \cref{fig:objective_function}.
It can also be seen that in most cases, a significant update of the Lagrange multipliers is accompanied by a substantial jump, negative or positive, in the value of the objective function.
This is signaled by the intersection data points between the changes in the norm level of the multipliers $\|\lambda_g\|$, the jumps of the objective, and the marked update steps.

\section{Algorithmic Scalability}\label{sec:scalability}
In this section we present weak scalability results for the 2d case presented in \cref{ss:2d}. 
These were carried out at \textit{HLRS} using the modern \textit{HPE Hawk} supercomputer. 
It consists of 5632 dual-socket nodes with the AMD EPYC 7742 processor. 
Each node has a total of 128 cores and 256GB of memory. 
The machine presents a 16-node hypercube connection topology, therefore the core counts are aimed at maximizing hypercube use, without significantly reducing the bandwidth. 
The grid partitioning is based on ParMETIS \cite{PARMETIS}.

\begin{figure}[!htbp]
	\centering
	\begin{subfigure}[t]{\textwidth}
		\centering
		\scalebox{0.46}{\footnotesize\includegraphics{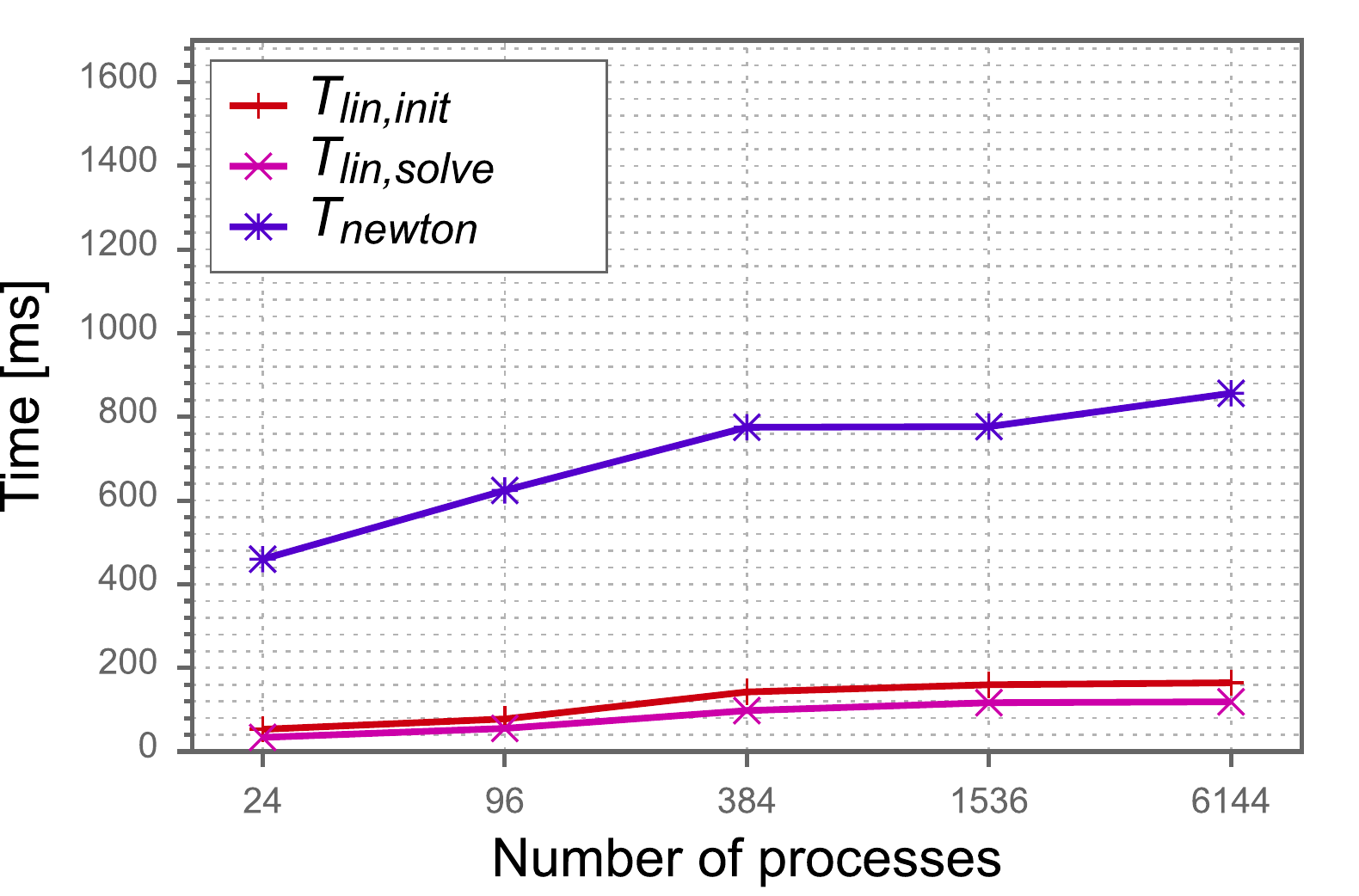}}
		\scalebox{0.46}{\footnotesize\includegraphics{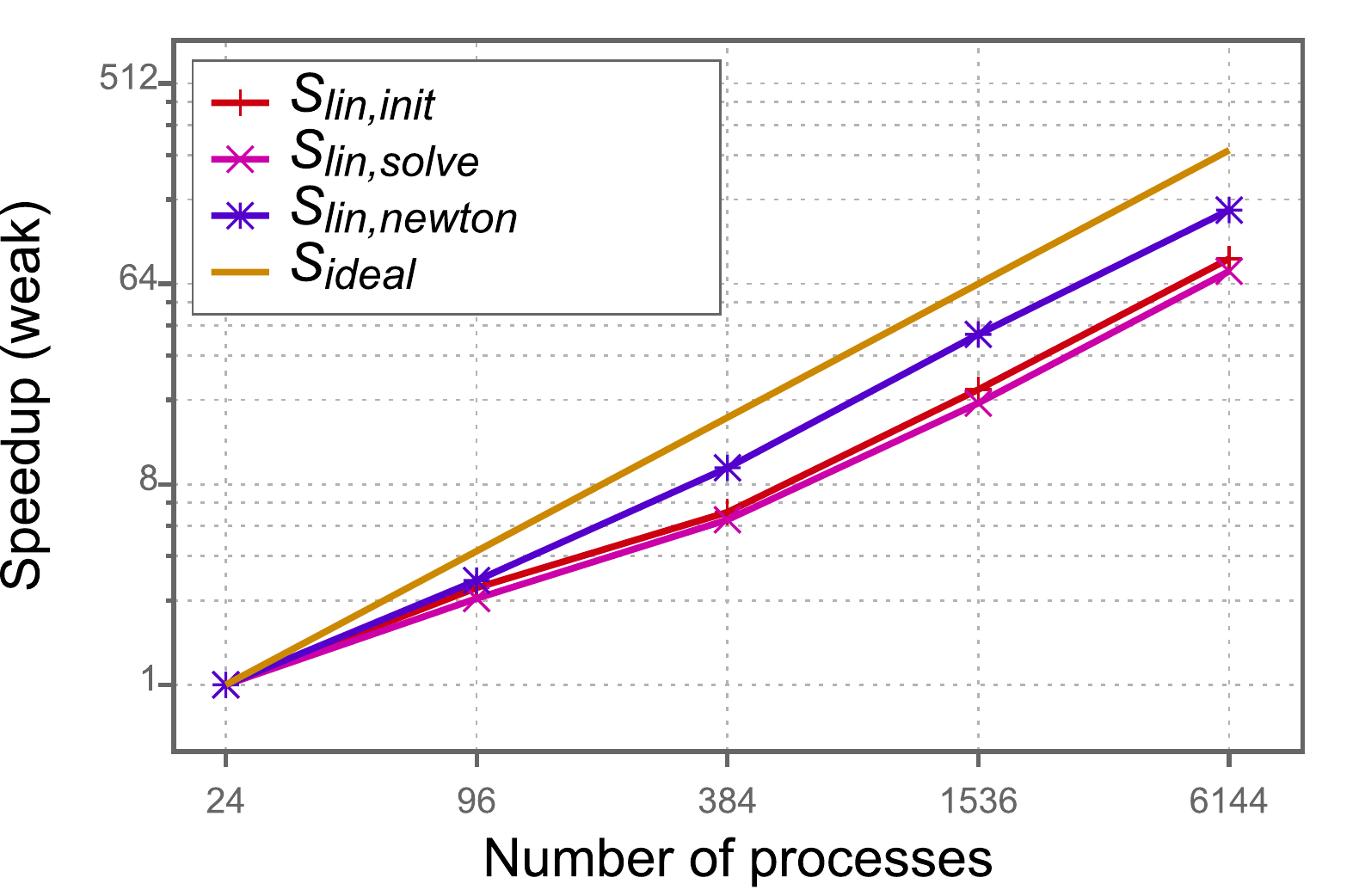}}		
		\caption{Extension equation nonlinear problem}
		\label{tab:weakDef}
	\end{subfigure}	
	\centering
	\begin{subfigure}[t]{\textwidth}
		\centering
		\scalebox{0.46}{\footnotesize\includegraphics{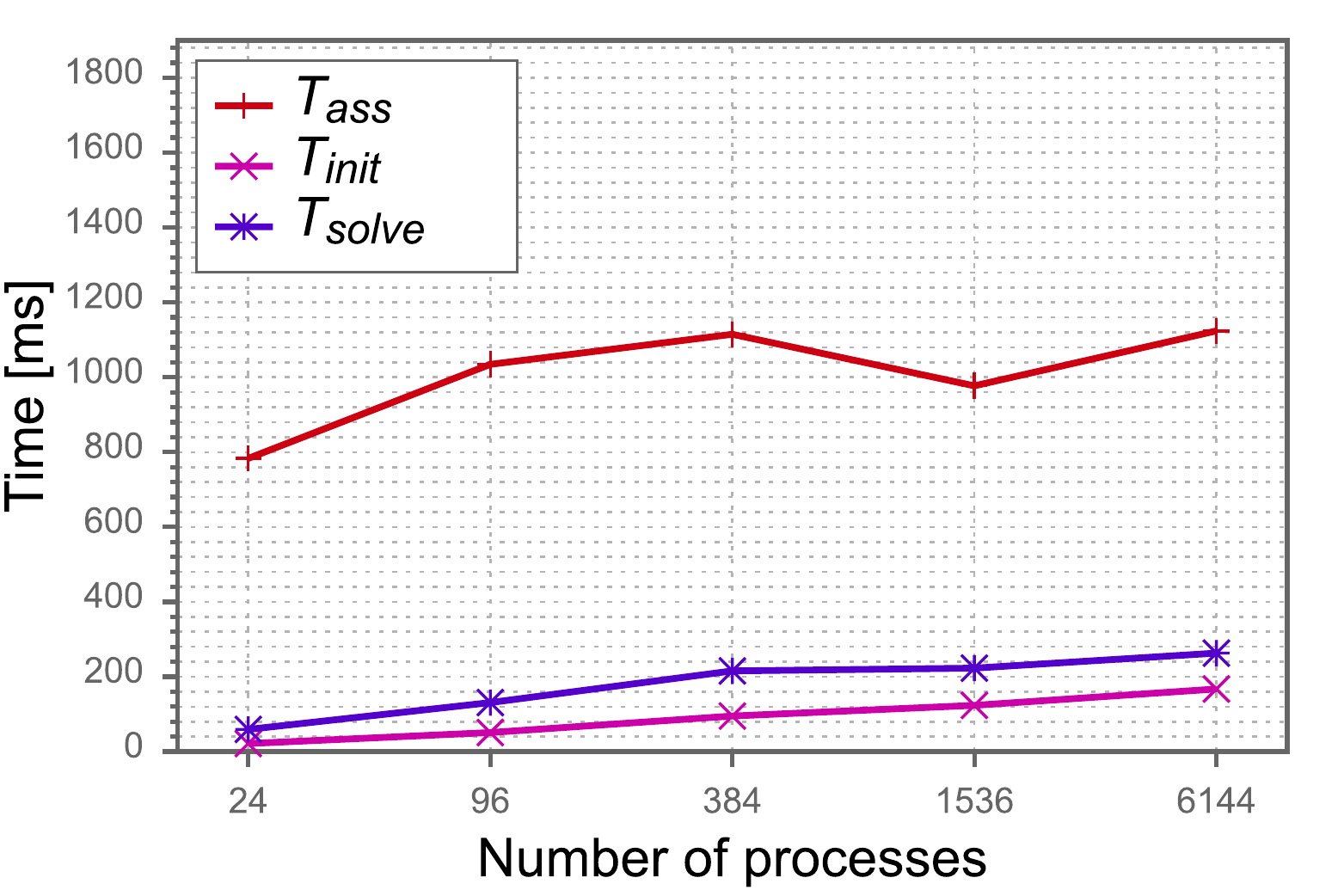}}
		\scalebox{0.46}{\footnotesize\includegraphics{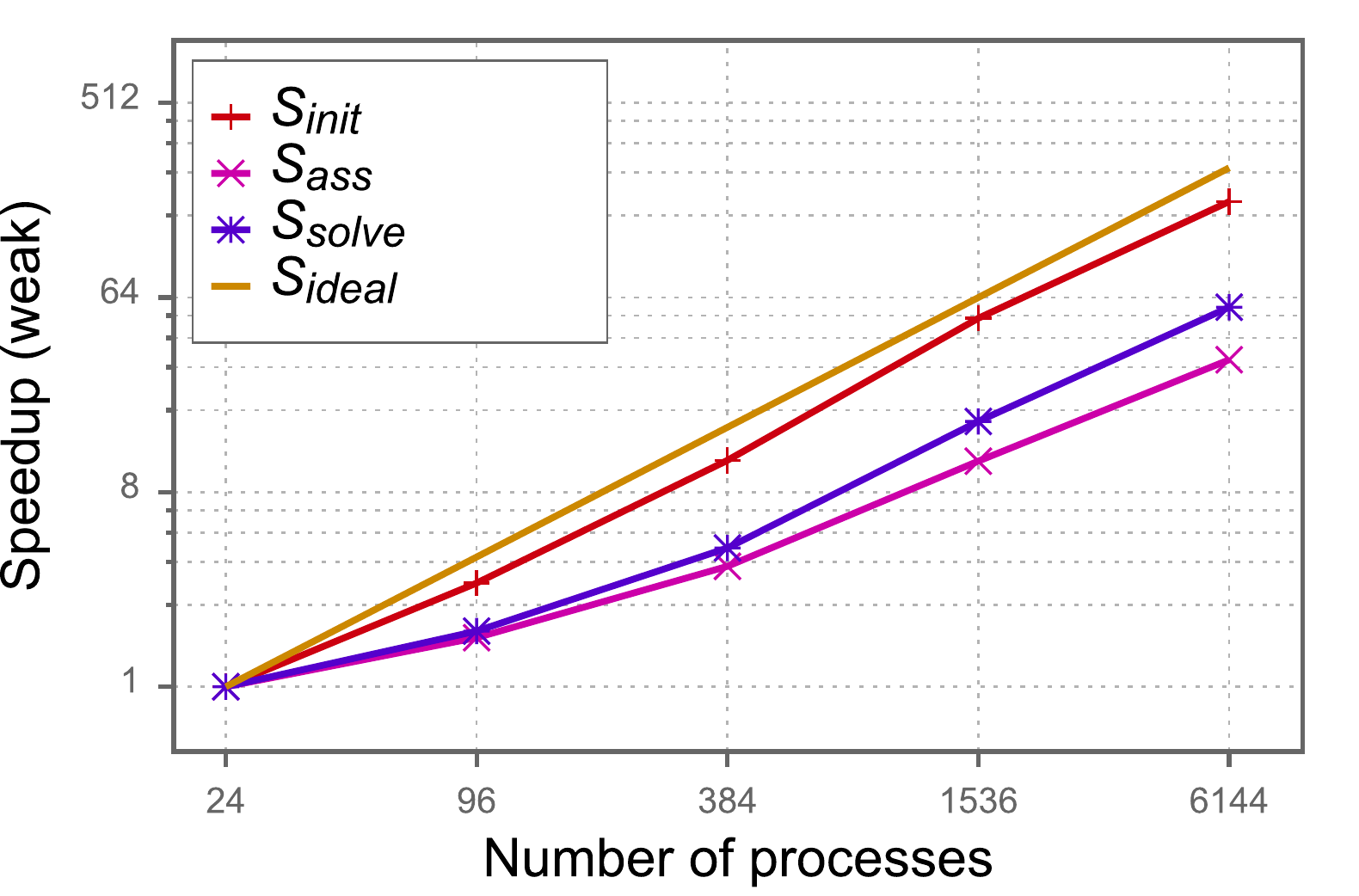}}		
		\caption{Adjoint displacement linear problem}
		\label{tab:weakAdj}
	\end{subfigure}
	\centering
	\begin{subfigure}[t]{\textwidth}
		\centering
		\scalebox{0.46}{\footnotesize\includegraphics{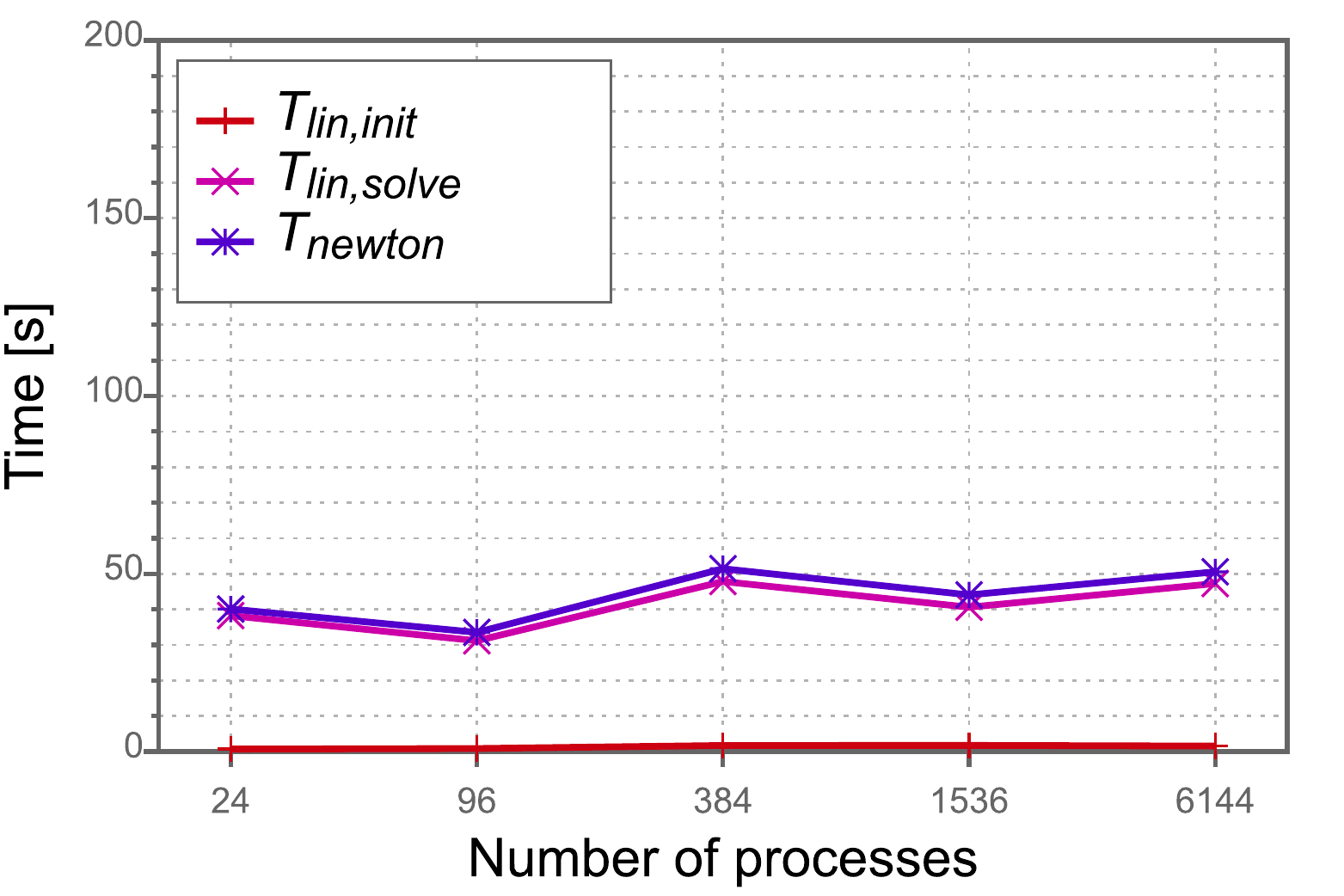}}
		\scalebox{0.46}{\footnotesize\includegraphics{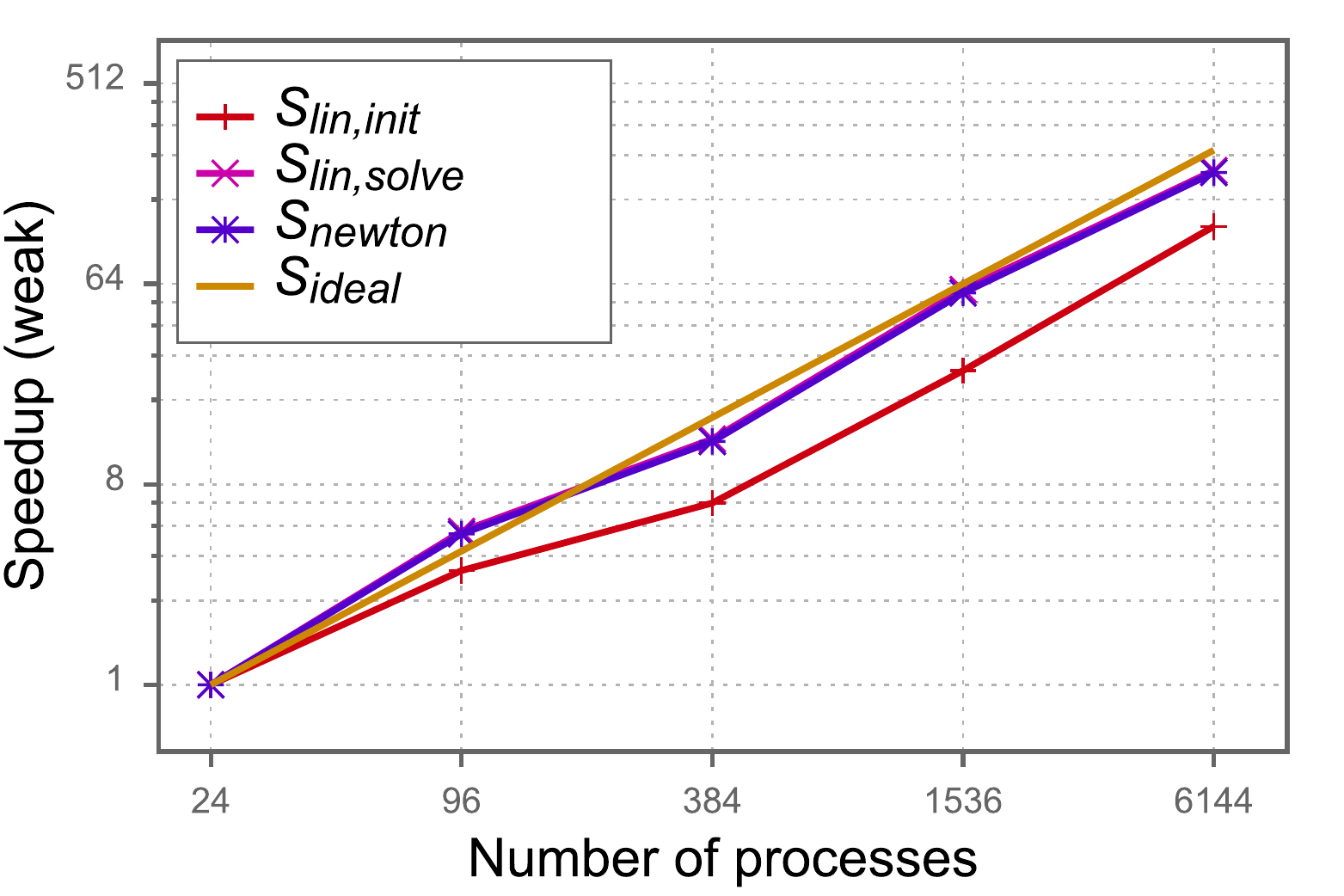}}	
		\caption{State equation nonlinear problem}
		\label{tab:weakGradExt}
	\end{subfigure}
	\vskip\baselineskip
	\begin{subfigure}[t]{\textwidth}
	\centering
	\begin{tabular}{rrrccc}
	\\
	\toprule
	Procs  & Refs  & NumElems  & Adjoint Displacement  & Displacement Field   & State Equation\\
	&       &        &  (linear solver)      & (linear/nonlinear solver)    & (linear/nonlinear solver)\\ 
	\midrule
	24 & 4 & 105472 & 56 & 9/12 & 238/12 \\
	96 & 5 & 421888 & 70 & 9/12 & 199/12 \\
	384 & 6 & 1687552 & 77 & 9/12 & 211/13 \\
	1536 & 7 & 6750208 & 77 & 9/12 & 194/13 \\
	6144 & 8 & 27000832 & 70 & 8/10 & 233/12 \\
	\bottomrule
	\end{tabular}
	\caption{Accumulated iteration counts for the weak scaling study}
	\label{tab:weakScaling}
	\end{subfigure}
	\caption{Weak Scaling: For the first three optimization steps, the accumulated wallclock time is shown for: (a) the nonlinear extension equation, (b) the derivative of the objective function with respect to the displacement field. In (c), the accumulated iteration counts are presented for the geometric multigrid preconditioned linear solver of the shape derivative, the number of Newton steps and linear iterations necessary to solve the extension equation and its linearization.}
	\label{fig:weak_scaling}
\end{figure}

\Cref{fig:weak_scaling} shows accumulated wallclock times, gained speedup relative to 24 cores and iteration counts for the first three optimization steps. 
These results are shown for the nonlinear extension operator \cref{eq:weak_deform_equation}, state \cref{eq:ns_1,eq:ns_2,eq:ns_3,eq:ns_4,eq:ns_5}, and the adjoint displacement \cref{eq:weak_adjoint_deform_equation}. 
A $P_1$ finite element discretization is used for extension operator and its adjoint equation, while mixed $P_1-P_1$ shape functions are used for the state equation Navier-Stokes equations. 
The nonlinear problems are solved using Newton's method altogether with a BiCGStab solver for the underlying linearizations.  
The linear solver is preconditioned with the geometric multigrid, which uses a V-cycle, 3 pre- and postsmoothing steps, and an LU base solver gathered on a single core. 
The error reduction is set to an absolute of $10^{-14}$ and relative $10^{-8}$ for the nonlinear solver, and $10^{-12}$ and $10^{-3}$ for the linear solvers. 
The purely linear problem has a relative and absolute reduction of $10^{-16}$. 
A Jacobi smoother is used within the geometric multigrid for the extension equation and the derivative, whereas the Navier-Stokes equation solver features an ILU smoother (see for instance \cite{wittum89}).
The results presented start at 24 cores, with a fourfold increase for each mesh refinement.
The studies show scalability and speedup for up to 6144 cores and more than 27 million triangular elements.
Given that mesh refinements are performed for each core count increase, a different geometrical problem is solved therefore differences in iteration counts are expected. 
However, even for a significant increase in the number of geometric elements the iteration counts for the linear problems remain within moderate bounds. Moreover, its important to point out that the total number of DOFs solved within the presented PDEs in \cref{fig:weak_scaling} increases from about 783k to 189 million. While the total number of DOFs solved in one optimization step is close to 300 million.

Together with the grid independence study for the outer optimization routine we thus obtain weak scalability of the overall method.

\section{Conclusion}
\label{sec:conclusion}
In this article we presented an optimization methodology which relies on the self-adaption of the extension operator within the method of mappings.
The results show that large deformations with respect to the reference configuration are possible while preserving mesh quality.
This has been studied in situations where singularities during the shape optimization process have to be smoothed out and newly generated.
It has been demonstrated that these two effects are particularly important to be tackled for applications of hierarchical multigrid solvers when experiments from fluid dynamics are considered.

The method's scalability and grid independency has been illustrated with the results of \cref{sec:scalability,ss:grid_stdy}. 
Grid independence is necessary for applications where a high level of refinement is needed, since it guarantees that the same optimal shape is obtained regardless of the number of elements. 
This becomes particularly important for the weak scalability, where the grid is refined on each core count increase. 
The results shown in \cref{fig:weak_scaling}, in combination with the ones of \cref{fig:gs_shapes,fig:gs_tips}, establish a proof of concept for industrial applicability, where a high number of DOFs are expected.
Overall, in this article we have presented an algorithm towards scalable shape optimization for large scale problems with the potential to work reliably also in complex geometric situations.

\subsubsection*{Acknowledgment}
\vspace{-0.2em}
{Computing time on the supercomputer {\it Hawk} at HLRS under the grant ShapeOptCompMat (ACID 44171, Shape Optimization for 3d Composite Material Models) is gratefully acknowledged.\\
	
\noindent
The current work is part of the research training group ‘Simulation-Based Design Optimization of Dynamic Systems Under Uncertainties’ (SENSUS) funded by the state of Hamburg under the aegis of the Landesforschungsförderungs-Project LFF-GK11.\\

\noindent
The authors acknowledge the support by the Deutsche Forschungsgemeinschaft (DFG) within the Research Training Group GRK 2583 “Modeling, Simulation and Optimization of Fluid Dynamic Applications”.}

\bibliographystyle{plain}      
\bibliography{nonlinext}   

\end{document}